\documentclass[12pt,leqno]{article}

\usepackage{amsmath,amssymb,amsthm}
\usepackage[all]{xy}

\makeatletter

\newtheorem*{theorem-main}{Theorem~\ref{main.thm}}

\newtheorem{theorem}{Theorem}[section]
\newtheorem{proposition}[theorem]{Proposition}
\newtheorem{lemma}[theorem]{Lemma}
\newtheorem{corollary}[theorem]{Corollary}

\theoremstyle{definition}

\newtheorem{example}[theorem]{Example}
\newtheorem{definition}[theorem]{Definition}

\theoremstyle{remark}

\theoremstyle{remark}

\def\({{\rm (}}
\def\){{\rm )}}

\let\Mathrm\operator@font
\let\Cal\mathcal
\let\Bbb\mathbb
\newcommand{\fm}{\ensuremath{\mathfrak m}}
\newcommand{\fn}{\ensuremath{\mathfrak n}}

\def\standop#1{\mathop{\Mathrm #1}\nolimits}
\def\difstop#1#2{\expandafter\def\csname #1\endcsname{\standop{#2}}}
\def\defstop#1{\difstop{#1}{#1}}

\defstop{AB}
\defstop{ann}
\defstop{Ass}
\defstop{Add}
\defstop{Alt}
\defstop{Ass}
\defstop{Assh}
\defstop{add}

\defstop{CMFI}
\defstop{codim}
\defstop{Coh}
\defstop{coht}
\defstop{Coker}
\defstop{Cone}
\defstop{Cl}
\defstop{Cox}
\defstop{cosk}
\defstop{cd}
\def\Cdim{\mathop{C\Mathrm{dim}}\nolimits}
\def\CCdim{\mathop{\C\Mathrm{dim}}\nolimits}

\defstop{depth}
\defstop{Der}
\defstop{Div}
\defstop{div}

\def\da{^\dagger}
\def\dda{^\ddagger}
\def\ddd{^{\dagger\ddagger}}
\def\dddd{^{\ddagger\dagger}}

\defstop{EM}
\defstop{embdim}
\defstop{End}
\defstop{ev}
\defstop{Ext}

\def\fg{_{\standop{fg}}}
\defstop{Flat}
\defstop{Func}
\defstop{Fpqc}

\defstop{GD}
\defstop{Good}
\defstop{Gal}

\difstop{height}{ht}
\defstop{Hom}
\def\hR{\hat{\Bbb R}}

\def\Id{\mathord{\Mathrm{Id}}}
\difstop{Image}{Im}
\defstop{ind}
\defstop{idim}

\defstop{Ker}

\defstop{Lch}
\defstop{length}
\defstop{Lin}
\defstop{Lqc}
\defstop{lqc}
\defstop{LQ}
\defstop{LM}
\defstop{Lie}
\defstop{Loc}
\def\Lop{_{\Lambda\op}}

\defstop{Mat}
\defstop{Max}
\defstop{Min}
\defstop{Mod}
\difstop{md}{mod}
\defstop{Mor}
\defstop{MCM}
\defstop{Map}
\difstop{mred}{red}

\defstop{Nerve}
\defstop{NonSFR}
\defstop{NonCMFI}
\defstop{NonNor}
\defstop{Nor}

\defstop{ob}
\defstop{Ob}
\def\op{^{\standop{op}}}

\defstop{PA}
\defstop{PM}
\defstop{PR}
\defstop{Proj}
\defstop{Prin}
\defstop{Pic}
\defstop{Push}

\defstop{Qch}
\defstop{qch}

\defstop{rad}
\defstop{rank}

\defstop{res}
\defstop{Reg}
\defstop{Ref}
\def\RHom{\standop{\mathbf{R}Hom}}

\defstop{Spec}
\defstop{supp}
\defstop{Supp}
\defstop{Sym}
\defstop{Sing}
\defstop{SFR}
\defstop{Soc}
\defstop{Sh}

\defstop{syz}
\defstop{Syz}
\defstop{SPH}
\def\sph#1{{}^{\perp_{#1}}}
\def\smd{\mathop{\text{\underline{$\Mathrm mod$}}}\nolimits}
\difstop{sm}{sum}

\difstop{tdeg}{trans.deg}

\defstop{Tor}
\difstop{trace}{tr}
\defstop{TF}
\defstop{tf}
\defstop{Tr}

\defstop{up}
\defstop{UP}

\defstop{Zar}

\def\an#1{^{\langle #1\rangle}}
\def\bra#1{^{[#1]}}

\def\A{\Cal A}
\def\B{\Cal B}
\def\C{\Cal C}
\def\E{\Cal E}

\def\I{\Cal I}

\def\M{\Cal M}
\def\N{\Cal N}
\def\O{\Cal O}

\def\Q{\Cal Q}

\def\V{\Cal V}
\def\X{\Cal X}
\def\Y{\Cal Y}
\def\Z{\Cal Z}

\def\fm{\mathfrak{m}}

\def\uHom{\mathop{\text{\underline{$\Mathrm Hom$}}}\nolimits}
\def\uEnd{\mathop{\text{\underline{$\Mathrm End$}}}\nolimits}

\def\sdarrow#1{\downarrow\hbox to 0pt{\scriptsize$#1$\hss}}
\def\suarrow#1{\uparrow\hbox to 0pt{\scriptsize$#1$\hss}}
\def\ssearrow#1{\searrow\hbox to 0pt{\scriptsize$#1$\hss}}

\def\section{\@startsection{section}{1}{\z@ }%
  {-3.5ex plus -1ex minus -.2ex}{2.3ex plus .2ex}{\bf }}

\long\def\refname{\par\kern -3ex
  \begin{center}\rm R\sc{eferences}\end{center}\par\kern 
  -2ex}

\def\@seccntformat#1{\csname the#1\endcsname.\quad}

\def\@@@sect#1#2#3#4#5#6[#7]#8{%
  \ifnum #2>\c@secnumdepth 
  \def \@svsec {}\else \refstepcounter {#1}%
  \def\@svsec{}
  \fi 
  \@tempskipa #5\relax 
  \ifdim \@tempskipa >\z@ 
  \begingroup #6\relax \@hangfrom {\hskip #3\relax 
    \@svsec}{\interlinepenalty \@M #8\par }\endgroup 
  \csname #1mark\endcsname {#7}
  \else 
  \def \@svsechd {#6\hskip #3\@svsec #8\csname #1mark\endcsname {#7}}
  \fi \@xsect {#5}}

\def\@@@startsection#1#2#3#4#5#6{%
  \if@noskipsec \leavevmode \fi \par \@tempskipa #4\relax \@afterindenttrue 
  \ifdim \@tempskipa <\z@ \@tempskipa -\@tempskipa \@afterindentfalse 
  \fi \if@nobreak \everypar {}\else \addpenalty {\@secpenalty }\addvspace 
  {\@tempskipa }\fi \@ifstar {\@ssect {#3}{#4}{#5}{#6}}{\@dblarg 
    {\@@@sect {#1}{#2}{#3}{#4}{#5}{#6}}}}

\def\theparagraph{\thesection.\arabic{paragraph}}
\def\aparagraph{\@@@startsection{paragraph}{2}{\z@ }%
  {1.75ex plus .2ex minus .15ex}{-1em}{\bf(\theparagraph) } }
\def\paragraph{\@@@startsection{paragraph}{2}{\z@ }%
  {1.75ex plus .2ex minus .15ex}{-1em}{}{\bf(\theparagraph)} }

\let\c@theorem\c@paragraph

\advance\textheight 10pt
\advance\textwidth 6pt
\title{Canonical and $n$-canonical modules on a Noetherian algebra}
\author{M{\sc itsuyasu} H{\sc ashimoto}}
\date{\normalsize
  Department of Mathematics, Okayama University\\
  Okayama 700--8530, JAPAN\\
  {\small \tt mh@okayama-u.ac.jp}\\
~\\
Dedicated to Professor Shiro Goto on the occasion of his 70th birthday
}

\begin{document}

\maketitle

\begin{abstract}
  We define canonical and $n$-canonical modules on a module-finite
  algebra over a Noether commutative ring and study their basic
  properties.
  Using $n$-canonical modules, we generalize a theorem on $(n,C)$-syzygy
  by Araya and Iima which generalize a well-known theorem on syzygies
  by Evans and Griffith.
  Among others, we prove a non-commutative version of Aoyama's theorem
  which states that a canonical module descends with respect to a flat
  local homomorphism.
  We also prove the codimension two-argument for modules over
  a coherent sheaf of algebras with a $2$-canonical module,
  generalizing a result of the author.
\end{abstract}

\tableofcontents

\section{Introduction}

\paragraph
In \cite{EG}, Evans and Griffith proved a criterion of a finite module
over a Noetherian commutative ring $R$ to be an $n$th syzygy.
This was generalized to a theorem on $(n,C)$-syzygy for a semidualizing module $C$
over $R$ by Araya and Iima \cite{AI}.
The main purpose of this paper is to prove a generalization of these results
in the following settings:
the ring $R$ is now a finite $R$-algebra $\Lambda$, which may not be commutative;
and $C$ is an $n$-canonical module.

\paragraph
The notion of $n$-canonical module was introduced in \cite{Hashimoto12} in an
algebro-geometric situation.
The criterion for a module to be an $n$th syzygy for $n=1,2$ by Evans--Griffith
was generalized using $n$-canonical modules there, and the standard
\lq codimension-two argument'
(see e.g., \cite[(1.12)]{Hartshorne4})
was also generalized to a theorem on
schemes with $2$-canonical modules \cite[(7.34)]{Hashimoto12}.
We also generalize this result to a theorem on modules over
non-commutative sheaves of algebras (Proposition~\ref{codim-two.prop}).

\paragraph
Let $(R,\fm)$ be a complete semilocal Noetherian ring, and $\Lambda\neq 0$ a module-finite
$R$-algebra.
Let $\Bbb I$ be a dualizing complex of $R$.
Then $\RHom_R(\Lambda,\Bbb I)$ is a dualizing complex of $\Lambda$.
Its lowest
non-vanishing cohomology is denoted by $K_\Lambda$, and is called the
canonical module of $\Lambda$.
If $(R,\fm)$ is semilocal but not complete, then
a $\Lambda$-bimodule is called a canonical module if it is the canonical module
after completion.
An $n$-canonical module is defined using the canonical module.
A finite right (resp.\ left, bi-)module $C$ of $\Lambda$ is said to be $n$-canonical
over $R$ 
if (1) $C$ satisfies Serre's $(S_n')$ condition as an $R$-module, that is,
for any $P\in\Spec R$, $\depth_{R_P} C_P\geq \min(n,\dim R_P)$.
(2) If $P\in\Supp_R C$ with $\dim R_P < n$, then $\widehat{C_P}$ is isomorphic to
$K_{\widehat{\Lambda_P}}$ as a right (left, bi-) module of $\widehat{\Lambda_P}$,
where $\widehat{\Lambda_P}$ is the $PR_P$-adic completion of $\Lambda_P$.

\paragraph
In order to study non-commutative $n$-canonical modules, we study some
non-commutative analogue of the theory of canonical modules developed by
Aoyama \cite{Aoyama}, Aoyama--Goto \cite{AG}, and Ogoma \cite{Ogoma}
in commutative algebra.
Among them, we prove an analogue of Aoyama's theorem \cite{Aoyama} which
states that the canonical module descends with respect to flat homomorphisms
(Theorem~\ref{Aoyama.thm}).

\paragraph
Our main theorem is the following.

\begin{theorem-main}[cf.~{\cite[(3.8)]{EG}, \cite[(3.1)]{AI}}]
  Let $R$ be a Noetherian commutative ring, and $\Lambda$ a module-finite
  $R$-algebra, which may not be commutative.
  Let $n\geq 1$, and $C$ be a right $n$-canonical $\Lambda$-module.
  Set $\Gamma=\End\Lop C$.
  Let $M\in\md C$.
  Then the following are equivalent.
  \begin{enumerate}
  \item[\bf 1] $M\in \TF(n,C)$.
  \item[\bf 2] $M\in \UP(n,C)$.
  \item[\bf 3] $M\in\Syz(n,C)$.
  \item[\bf 4] $M\in (S_n')_C$.
  \end{enumerate}
\end{theorem-main}

Here $M\in (S_n')_C$ means that $\Supp_R M\subset \Supp_R C$, and
for any $P\in\Spec R$, $\depth M_P\geq \min(n,\dim R_P)$, and this is a
(modified) Serre's condition.
$M\in\Syz(n,C)$ means $M$ is an $(n,C)$-syzygy.
$M\in\UP(n,C)$ means existence of an exact sequence
\[
0\rightarrow M\rightarrow C^0\rightarrow C^1\rightarrow\cdots\rightarrow C^{n-1}
\]
which is still exact after applying $(?)^\dagger=\Hom_{\Lop}(?,C)$.

\paragraph
The condition $M\in\TF(n,C)$ is a modified version of Takahashi's condition
``$M$ is $n$-$C$-torsion free'' \cite{Takahashi}.
Under the assumptions of the theorem, let $(?)^\dagger=\Hom_{\Lop}(?,C)$,
$\Gamma=\End_{\Lop} C$, and $(?)^\ddagger=\Hom_\Gamma(?,C)$.
We say that $M\in\TF(1,C)$ (resp.\ $M\in\TF(2,C)$)
if the canonical map $\lambda_M:M\rightarrow M^{\dagger\ddagger}$ is
injective (resp.\ bijective).
If $n\geq 3$, we say that $M\in\TF(n,C)$ if $M\in\TF(2,C)$, and
$\Ext^i_\Gamma(M^\dagger,C)=0$ for $1\leq i\leq n-2$, see Definition~\ref{(n,C)-TF.def}.
Even if $\Lambda$ is a commutative ring, a non-commutative ring
$\Gamma$ appears in a natural way, so even in this case, the definition
is slightly different from Takahashi's original one.
We prove that $\TF(n,C)=\UP(n,C)$ in general (Lemma~\ref{Takahashi.lem}).
This is a modified version of Takahashi's result \cite[(3.2)]{Takahashi}.

\paragraph
As an application of the main theorem,
we formulate and prove a different form of the existence of
$n$-$C$-spherical approximations by 
Takahashi \cite{Takahashi}, using $n$-canonical modules,
see Corollary~\ref{main.cor} and Corollary~\ref{main2.cor}.
Our results are not strong enough to deduce \cite[Corollary~5.8]{Takahashi}
in commutative case.
For related categorical results, see below.

\paragraph
Section~\ref{quotients.sec} is preliminaries on the depth and Serre's conditions on
modules.
In Section~\ref{spherical.sec}, we discuss $\Cal X_{n,m}$-approximation,
which is a categorical abstraction of approximations of
modules appeared in \cite{Takahashi}.
Everything is done categorically here, and Theorem~\ref{X_{m,n}-approximation.thm} is
an abstraction of 
\cite[(3.5)]{Takahashi}, in view of the fact that 
$\TF(n,C)=\UP(n,C)$ in general (Lemma~\ref{Takahashi.lem}).
In Section~\ref{(n,C)-TF.sec}, we discuss $\TF(n,C)$, and
prove Lemma~\ref{Takahashi.lem} and related lemmas.
In Section~\ref{canonical.sec}, we define the canonical module of a module-finite
algebra $\Lambda$
over a Noetherian commutative ring $R$, and prove some basic properties.
In Section~\ref{n-canonical.sec}, we define the $n$-canonical module of $\Lambda$, and
prove some basic properties, generalizing some constructions and results in
\cite[Section~7]{Hashimoto12}.
In Section~\ref{Aoyama.sec}, we prove a non-commutative version of Aoayama's theorem
which says that the canonical module descends with respect to flat local homomorphisms
(Theorem~\ref{Aoyama.thm}).
As a corollary, as in the commutative case, we immediately have that a localization
of a canonical module is again a canonical module.
This is important in Section~\ref{Evans-Griffith.sec}.
In Section~\ref{Evans-Griffith.sec}, we prove Theorem~\ref{main.thm}, and
the related
results on $n$-$C$-spherical approximations
(Corollary~\ref{main.cor}, Corollary~\ref{main2.cor}) as its corollaries.
Before these, we prove non-commutative analogues of the theorems of Schenzel and
Aoyama--Goto \cite[(2.2), (2.3)]{AG} on the Cohen--Macaulayness of the canonical module
(Proposition~\ref{AG.prop} and Corollary~\ref{AG.cor}).
In section~\ref{Frobenius.sec}, we define and discuss non-commutative, higher-dimensional
symmetric, Frobenius, and quasi-Frobenius algebras and their non-Cohen--Macaulay versions.
In commutative algebra, the non-Cohen--Macaulay version of Gorenstein ring is known as
quasi-Gorenstein rings.
What we discuss here is a non-commutative version of such rings.
Scheja and Storch \cite{SS} discussed a relative notion, and our definition is
absolute in the sense that it is independent of the choice of $R$.
If $R$ is local, our quasi-Frobenius property agrees with
Gorensteinness discussed by Goto and Nishida \cite{GN},
see Proposition~\ref{GN.prop} and Corollary~\ref{GN.cor}.
In Section~\ref{codim-two.sec}, we show that the codimension-two argument using the
existence of $2$-canonical modules in \cite{Hashimoto12} is still valid in
non-commutative settings.

\paragraph
Acknowledgments:
Special thanks are due to Professor
Osamu Iyama for valuable advice and
discussion.
Special thanks are also due to
Professor Tokuji Araya.
This work was motivated by his advice, and
Proposition~\ref{AG.prop}
is an outcome of discussion with him.

The author is also grateful to
Professor Kei-ichiro Iima,
Professor Takesi Kawasaki,
Professor Ryo Takahashi,
Professor Kohji Yanagawa,
and
Professor Yuji Yoshino
for valuable advice.

\section{Preliminaries}\label{quotients.sec}

\paragraph
Unless otherwise specified, a module means a left module.
Let $B$ be a ring.
$\Hom_B$ or $\Ext_B$ mean the $\Hom$ or $\Ext$ for left $B$-modules.
$B\op$ denotes the opposite ring of $B$, so a $B\op$-module
is nothing but a right $B$-module.
Let $B\Mod$ denote the category of $B$-modules.
$B\op\Mod$ is also denoted by $\Mod B$.
For a left (resp.\ right)
Noetherian ring $B$, $B\md$ (resp.\ $\md B$)
denotes the full subcategory
of $B\Mod$ (resp.\ $\Mod B$)
consisting of finitely generated left (resp.\ right) $B$-modules.

\paragraph
For derived categories, we employ standard notation found in
\cite{Hartshorne2}.

For an abelian category $\A$,
$D(\A)$ denotes the unbounded derived category of $\A$.
For a plump subcategory
(that is, a full subcategory which is closed under kernels, cokernels, and
extensions) $\B$ of $\A$, 
$D_{\B}(\A)$ denotes the triangulated subcategory of $D(\A)$ consisting
of objects $\Bbb F$ such that $H^i(\Bbb F)\in \B$ for any $i$.
For a ring $B$, 
We denote $D(B\Mod)$ by $D(B)$, and 
$D_{B\md}(B\Mod)$ by $D\fg(B)$ (if $B$ is left Noetherian).

\paragraph
Throughout the paper, let $R$ denote a commutative Noetherian ring.
If $R$ is semilocal (resp.\ local) and $\fm$ its Jacobson radical,
then we say that $(R,\fm)$ is semilocal (resp.\ local).
We say that $(R,\fm,k)$ is semilocal (resp.\ local) if $(R,\fm)$ is
semilocal (resp.\ local) and $k=R/\fm$.

\paragraph
We set $\hR:=\Bbb R\cup\{\infty,-\infty\}$ and
consider that $-\infty<\Bbb R<\infty$.
As a convention, 
for a subset $\Gamma$ of $\hR$,
$\inf\Gamma$ means $\inf (\Gamma\cup\{\infty\})$, which exists uniquely
as an element of $\hR$.
Similarly for $\sup$.

\paragraph
For an ideal $I$ of $R$ and $M\in\md R$, we define
\[
\depth_R(I,M):=\inf\{i\in\Bbb Z\mid \Ext^i_R(R/I,M)\neq 0\},
\]
and call it the $I$-depth of $M$ \cite[section~16]{CRT}.
It is also called the $M$-grade of $I$ \cite[(6.2.4)]{BS}.
When $(R,\fm)$ is semilocal, we 
denote $\depth(\fm,M)$ by $\depth_R M$ or $\depth M$,
and call it the depth of $M$.

\begin{lemma}\label{depth(I,M).lem}
  The following functions on $M$ \(with valued in $\hat R$\)
  are equal for an ideal $I$ of $R$.
  \begin{enumerate}
  \item[\bf 1] $\depth_R (I,M)$;
  \item[\bf 2] $\inf_{P\in V(I)}\depth_{R_P} M_P$,
    where $V(I)=\{P\in\Spec R\mid P\supset I\}$;
  \item[\bf 3] $\inf\{i\in\Bbb Z \mid H^i_I(M)\neq 0\}$;
  \item[\bf 4] $\infty$ if $M=IM$, and otherwise,
    the length of any maximal $M$-sequence in $I$.
  \item[\bf 5] Any function $\phi$ such that
    \begin{enumerate}
    \item[\bf a] $\phi(M)=\infty$ if $M=IM$.
    \item[\bf b]  $\phi(M)=0$ if $\Hom_R(R/I,M)\neq 0$.
    \item[\bf c] $\phi(M)=\phi(M/aM)+1$ if $a\in I$ is a
      nonzerodivisor on $M$.
  \end{enumerate}
  \end{enumerate}
\end{lemma}

  \begin{proof}
    We omit the proof, and refer the reader to \cite[section~16]{CRT},
    \cite[(6.2.7)]{BS}.
  \end{proof}

\paragraph
For a subset $F$ of $X=\Spec R$, we define $\codim F=\codim_XF$,
the {\em codimension}
of $F$ in $X$, by $\inf\{\height P\mid P\in F\}$.
So $\height I=\codim V(I)$ for an ideal $I$ of $R$.
For $M\in\md R$, we define $\codim M:=\codim\Supp_R M=\height \ann M$,
where $\ann$ denotes the annihilator.
For $n\geq 0$, we denote the set 
$\height^{-1}(n)=\{P\in \Spec R\mid \height P=n\}$ by $R\an n$.
For a subset $\Gamma$ of $\Bbb Z$, $R\an\Gamma$ means $
\height^{-1}(\Gamma)=\bigcup_{n\in \Gamma} R\an n$.
Moreover, we use notation such as $R\an{\leq 3}$, which stands for
$R\an{\{n\in\Bbb Z\mid n\leq 3\}}$.
For $M\in\md R$, the set of minimal primes of $M$ is denoted by $\Min M$.

We define $M\bra n:=\{P\in \Spec R\mid \depth M_P=n\}$.
Similarly, we use notation such as
$M\bra{<n}(=\{P\in \Spec R \mid \depth M_P<n\})$.

\paragraph
Let $M,N\in\md R$.
We say that $M$ satisfies the $(S_n^{N})^R$-condition or
$(S_n^{N})$-condition
if for any $P\in \Spec R$, $\depth_{R_P}M_P \geq\min(n,\dim_{R_P}N_P)$.
The $(S_n^R)^R$-condition or $(S_n^R)$-condition
is simply denoted by $(S_n')^R$ or $(S_n')$.
We say that $M$ satisfies the $(S_n)^R$-condition or $(S_n)$-condition if
$M$ satisfies the $(S_n^{M})$-condition.
$(S_n)$ (resp.\ $(S_n')$) is equivalent to say that for any
$P\in M\bra{<n}$,
$M_P$ is a Cohen--Macaulay (resp.\ maximal Cohen--Macaulay) $R_P$-module.
That is, $\depth M_P=\dim M_P$ (resp.\ $\depth M_P=\dim R_P$).
We consider that $(S_n^N)^R$ is a class of modules, and
also write $M\in(S_n^N)^R$ (or $M\in (S_n^N)$).

\begin{lemma}\label{depth-lemma.lem}
  Let $0\rightarrow L\rightarrow M\rightarrow N\rightarrow 0$
  be an exact sequence in $\md R$, and $n\geq 1$.
  \begin{enumerate}
  \item[\bf 1] If $L,N\in (S_n')$, then $M\in (S_n')$.
  \item[\bf 2] If $N\in (S_{n-1}')$ and $M\in (S_n')$, then $L\in (S_n')$.
  \end{enumerate}
\end{lemma}

\begin{proof}
  {\bf 1} follows from the depth lemma:
  \[
  \forall P\quad
  \depth_{R_P}M_P\geq
  \min(\depth_{R_P}L_P,\depth_{R_P}N_P),
  \]
  and the fact that maximal Cohen--Macaulay
  modules are closed under extensions.
  {\bf 2} is similar.
\end{proof}

\begin{corollary}\label{S_n-increase.cor}
  Let
  \[
  0\rightarrow M\rightarrow L_{n}\rightarrow\cdots \rightarrow
  L_1
  \]
  be an exact sequence in $\md R$, and assume that
  $L_i\in (S_i')$ for $1\leq i\leq n$.
  Then $M\in (S_n')$.
\end{corollary}

\begin{proof}
  This is proved using a repeated use of Lemma~\ref{depth-lemma.lem},
  {\bf 2}.
\end{proof}

\begin{lemma}[Acyclicity Lemma, {\cite[(1.8)]{PS}}]\label{gw.lem}
  Let $(R,\fm)$ be a Noetherian local ring, and
  \begin{equation}\label{L.eq}
  \Bbb L:
  0\rightarrow L_s\xrightarrow{\partial_s}
  L_{s-1}\xrightarrow{\partial_{s-1}}\rightarrow\cdots\rightarrow
  L_1\xrightarrow{\partial_1}L_0
  \end{equation}
  be a complex of $\md R$ such that
  \begin{enumerate}
  \item[\bf 1] For each $i\in\Bbb Z$ with $1\leq i\leq s$,
    $\depth L_i\geq i$.
  \item[\bf 2] For each $i\in\Bbb Z$ with $1\leq i\leq s$,
    $H_i(\Bbb L)\neq 0$ implies that $\depth H_i(\Bbb L)=0$.
  \end{enumerate}
  Then $\Bbb L$ is acyclic \(that is, $H_i(\Bbb L)=0$ for $i>0$\).
  \qed
\end{lemma}

\begin{lemma}[cf.~{\cite[(3.4)]{IW}}]\label{s_n-acyclicity.lem}
  Let {\rm(\ref{L.eq})} be a complex in $\md R$ such that
  \begin{enumerate}
  \item[\bf 1] For each $i\in\Bbb Z$ with $1\leq i\leq s$,
    $L_i\in (S_i')$.
  \item[\bf 2] For each $i\in\Bbb Z$ with $1\leq i\leq s$,
    $\codim H_i(\Bbb L)\geq s-i+1$.
  \end{enumerate}
  Then $\Bbb L$ is acyclic.
\end{lemma}

\begin{proof}
Using induction on $s$, we may assume that $H_i(\Bbb L)=0$ for $i>1$.
Assume that $\Bbb L$ is not acyclic.
Then $H_1(\Bbb L)\neq 0$, and we can take $P\in \Ass_R H_1(\Bbb L)$.
By assumption, $\height P \geq s$.
Now localize at $P$ and considering the complex
$\Bbb L_P$ over $R_P$, we get a contradiction by Lemma~\ref{gw.lem}.
\end{proof}

\begin{example}\label{injective-bijective.ex}
  Let $f:M\rightarrow N$ be a map in $\md R$.
  \begin{enumerate}
  \item[\bf 1] If $M\in (S_1')$ and
    $f_P$ is injective for $P\in R\an 0$, then $f$ is injective.
    Indeed, consider the complex
\[
0\rightarrow M\xrightarrow{f} N=L_0
\]
and apply
Lemma~\ref{s_n-acyclicity.lem}.
\item[\bf 2] (\cite[(5.11)]{LW}) If $M\in (S_2')$, $N\in (S_1')$, and
  $f_P$ is bijective for $P\in R\an{\leq 1}$, then $f$ is bijective.
  Consider the complex
  \[
  0\rightarrow M\xrightarrow {f} N\rightarrow 0=L_0
  \]
  this time.
  \end{enumerate}
\end{example}

\begin{lemma}\label{Ischebeck.lem}
  Let $(R,\fm)$ be a Noetherian local ring, and $N\in (S_n)^R$.
  If $P\in\Min N$ with $\dim R/P<n$, then we have
  \[
  \dim R/P = \depth N = \dim N<n.
  \]
  If, moreover, $N\in (S_n')^R$, then $\depth N=\dim R$.
\end{lemma}

\begin{proof}
  Ischebeck proved that if $M,N\in \md R$ and $i<\depth N-\dim M$,
  then $\Ext^i_R(M,N)=0$ \cite[(17.1)]{CRT}.
  As $\Ext^0_R(R/P,N)\neq 0$, we have that
  $\depth_R N\leq \dim R/P <n$.
  The rest is easy.
\end{proof}

\begin{corollary}\label{S_n.cor}
  Let $M\in (S_n)^R$ and $N\in (S_n')^R$.
  If $\Min M\subset \Min N$, then $M\in (S_n')^R$.
\end{corollary}

\begin{proof}
  Let $P\in M\bra{<n}$.
  As $M\in (S_n)$, $\depth M_P=\dim M_P$.
  Take $Q\in \Min M$ such that $Q\subset P$ and $\dim R_P/QR_P=\dim M_P < n$.
  As $\Min M\subset\Min N$, we have that $QR_P\in \Min N_P$.
  By Lemma~\ref{Ischebeck.lem},
  $\dim R_P=\dim R_P/QR_P=\depth M_P$, and hence $M\in (S_n')$.
\end{proof}

\begin{corollary}\label{(S_n)cap(S_1').cor}
  Let $n\geq 1$, and $R\in (S_n)$.
  Then for $M\in \md R$, we have that $(S_n')^R= (S_n)^R\cap (S_1')$.
\end{corollary}

\begin{proof}
  Obviously, $(S_n')^R\subset (S_n)^R \cap (S_1')$.
  For the converse, apply Corollary~\ref{S_n.cor} for $N=R$.
\end{proof}

\paragraph\label{(S_n')_N.par}
Let $M,N\in \md R$.
We say that
$M$ satisfies the $(S_n')_N$-condition, or $M\in (S_n')_N=(S_n')_N^R$,
if $M\in (S_n')$ and $\Supp_R M\subset \Supp_R N$.

\begin{lemma}\label{1-2.lem}
  Let $n\geq 1$, $N\in (S_n')$, and $M\in\md R$.
  Then the following are equivalent.
  \begin{enumerate}
    \item[\bf 1] $M\in (S_n')_N$.
    \item[\bf 2] $M\in (S_n)$ and $\Min M\subset \Min N$.
  \end{enumerate}
\end{lemma}

\begin{proof}
  {\bf 1$\Rightarrow$2}.
  As $(S_n')\subset (S_n)$, $M\in (S_n)$.
  As $M\in (S_n')$ with $n\geq 1$, $\Min M\subset \Min R$.
  By assumption, $\Min M\subset \Supp N$.
  So $\Min M\subset \Min R\cap \Supp N\subset \Min N$.

  {\bf 2$\Rightarrow$1}.
  $M\in (S_n')$ by Corollary~\ref{S_n.cor}.
  $\Supp M\subset \Supp N$ follows from $\Min M\subset \Min N$.
\end{proof}

\paragraph
There is another case that $(S_n)$ implies $(S_n')$.
An $R$-module $N$ is said to be {\em full} if $\Supp_R N=\Spec R$.
A finitely generated faithful $R$-module is full.

\begin{lemma}\label{S_n-S_n'.lem}
  Let $M,N\in\md R$.
  If $N$ is a full $R$-module,
  then $M$ satisfies $(S_n')$ condition if and only if
  $M$ satisfies $(S_n^N)$ condition.
  If $\ann_R N\subset \ann_R M$,
  then $M$ satisfies the $(S_n^N)^R$ condition if
  and only if $M$ satisfies the $(S_n')^{R/\ann_R N}$ condition.
\end{lemma}

\begin{proof}
  The first assertion is because
  $\dim N_P=\dim R_P$ for any $P\in\Spec R$.
  The second assertion follows from the first, because
  for an $R/\ann_R N$-module, $(S_n^N)^R$ and
  $(S_n^N)^{R/\ann_R N}$
  are the same thing.
\end{proof}
  
\begin{lemma}\label{depth_R=depth_S.lem}
  Let $I$ be an ideal of $R$, and $S$ a module-finite commutative
  $R$-algebra.
  For $M\in\md S$, we have that $\depth_R(I,M)=\depth_S(IS,M)$.
  In particular, if $R$ is semilocal, then $\depth_R M=\depth_S M$.
\end{lemma}
  
\begin{proof}
Note that $H^i_I(M)\cong H^i_{IS}(M)$ by \cite[(4.2.1)]{BS}.
By Lemma~\ref{depth(I,M).lem}, we get the lemma immediately.
\end{proof}

\begin{lemma}\label{finite-s_n.lem}
  Let $\varphi:R\rightarrow S$ be a finite homomorphism of rings,
  $M\in\md S$, and $n\geq 0$.
  \begin{enumerate}
  \item[\bf 1] If $M\in (S_n')^R$, then $M\in (S_n')^S$.
  \item[\bf 2]
    Assume that for any $Q\in\Min S$, $\varphi^{-1}(Q)\in\Min R$
    \(e.g., $S\in (S_1')^R$\).
    If $M\in (S_n')^S$, and
    $R_P$ is quasi-unmixed for any
    $P \in R\bra{<n}$, then $M\in (S_n')^R$.
  \end{enumerate}
\end{lemma}

\begin{proof}
  {\bf 1}.
  Let $Q\in M\bra{<n}$.
  Then $\depth_{R_P}M_P=\depth_{S_P}M_P\leq \depth_{S_Q}M_Q<n$ by
  Lemma~\ref{depth_R=depth_S.lem} and Lemma~\ref{depth(I,M).lem},
  where $P=\varphi^{-1}(Q)$.
  So $M_P$ is a maximal Cohen--Macaulay $R_P$-module by the $(S_n')_R$-property,
  and hence $\height Q\leq \height P=\depth_{R_P}M_P\leq \depth_{S_Q}M_Q$,
  and hence $M_Q$ is a maximal Cohen--Macaulay
  $S_Q$-module, and $M\in (S_n')_S$.

  {\bf 2}.
  Let $P\in\Spec R$, and $\depth_{R_P}M_P<n$.
  Then by
  Lemma~\ref{depth_R=depth_S.lem} and Lemma~\ref{depth(I,M).lem},
  there exists some $Q\in\Spec S$ such that
  $\varphi^{-1}(Q)=P$ and
  \[
  \depth_{S_Q}M_Q
  =
  \inf_{\varphi^{-1}(Q')=P}\depth_{S_{Q'}}M_{Q'}
  =
  \depth_{S_P}M_P = \depth_{R_P}M_P<n.
  \]
  Then $\height Q=\depth R_P M_P$.
  So it suffices to show $\height P=\height Q$.
  By assumption, $R_P$ is quasi-unmixed.
  So $R_P$ is equi-dimensional and universally catenary \cite[(31.6)]{CRT}.
  By \cite[(13.3.6)]{EGA-IV-3}, $\height P=\height Q$, as desired.
\end{proof}
  
\paragraph
We say that $R$ satisfies $(R_n)$ (resp.\ $(T_n)$) if
$R_P$ is regular (resp.\ Gorenstein) for $P\in R\an{\leq n}$.

\begin{lemma}\label{flat-S_n.lem}
  Let $\varphi:R\rightarrow S$ be a flat morphism between Noetherian
  rings, and $M\in\md R$.
  \begin{enumerate}
  \item If $M\in (S_n')^R$ and the ring $S_P/PS_P$ satisfies $(S_n)$ for
    $P\in\Spec R$, then $S\otimes_R M\in (S_n')^S$.
  \item If $\varphi$ is faithfully flat and $S\otimes_R M\in (S_n')^S$,
    then $M\in (S_n')^R$.
  \item If $R$ satisfies $(S_n)$ \(resp.\ $(T_n)$, $(R_n)$\) and
    $S_P/PS_P$ satisfies $(S_n)$ \(resp.\ $(T_n)$, $(R_n)$\)
    for $P\in\Spec R$,
    then $S$ satisfies $(S_n)$ \(resp.\ $(T_n)$, $(R_n)$\).
  \end{enumerate}
\end{lemma}

\begin{proof}
 Left to the reader (see \cite[(23.9)]{CRT}).
\end{proof}

\section{$\X_{n,m}$-approximation}\label{spherical.sec}

\paragraph
Let $\A$ be an abelian category, and $\C$ its additive subcategory
closed under direct summands.
Let $n\geq 0$.
We define
\[
\sph n\C:=\{a\in\A\mid\Ext^i_\A(a,c)=0\quad 1\leq i\leq n\}.
\]
Let $a\in\A$.
A sequence
\begin{equation}\label{n-spherical-cores.eq}
  \Bbb C:
  0\rightarrow a\rightarrow c^0\rightarrow c^1\rightarrow\cdots\rightarrow
c^{n-1}
\end{equation}
is said to be an $(n,\C)$-pushforward if it is exact
with $c^i\in\C$.
If in addition,
\[
\Bbb C\da:
0\leftarrow a\da\leftarrow (c^0)\da\leftarrow (c^1)\da\leftarrow
\cdots\leftarrow (c^{n-1})\da
\]
is exact for any $c\in \C$, where $(?)\da=\Hom_{\A}(?,c)$, we say that
$\Bbb C$ is a universal $(n,\C)$-pushforward.

If $a\in\A$ has an $(n,\C)$-pushforward, we say that
$a$ is an $(n,\C)$-syzygy, and we write $a\in \Syz(n,\C)$.
If $a\in \A$ has a universal $(n,\C)$-pushforward, we say that
$a\in\UP_{\A}(n,\C)=\UP(n,\C)$.
Obviously, $\UP_{\A}(n,\C)\subset\Syz_{\A}(n,\C)$.

\paragraph
We write $\X_{n,m}(\C)=\X_{n,m}:=\sph n\C\cap \UP(m,\C)$ for $n,m\geq 0$.
Also, for $a\neq 0$, we define
\begin{multline*}
\CCdim a=\inf\{m\in\Bbb Z_{\geq 0}\mid
\text{there is a resolution }\\
0\rightarrow c_m\rightarrow c_{m-1}\rightarrow
\cdots\rightarrow c_0\rightarrow a\rightarrow 0\}.
\end{multline*}
We define $\CCdim 0=-\infty$.
We define $\Y_n(\C)=\Y_n:=\{a\in\A\mid \CCdim a < n\}$.
A sequence $\Bbb E$ is said to be $\C$-exact if it is exact, and
$\A(\Bbb E,c)$ is also exact for each $c\in \C$.
Letting a $\C$-exact sequence an exact sequence, $\A$ is an exact
category, which we denote by $\A_{\C}$ in order to distinguish it from
the abelian category $\A$ (with the usual exact sequences).

\paragraph
Let $\C_0\subset\A$ be a subset.
Then $\sph n \C_0$, $\UP(n,\C_0)$, $\Cal X_{n,m}(\C_0)$,
$\mathop{\C_0\Mathrm{dim}}$, and $\Cal Y_n(\C_0)=\Cal Y_n$ mean
$\sph n \C$, $\UP(n,\C)$, $\Cal X_{n,m}(\C)$, $\CCdim$, and $\Y_n(\C)$,
respectively, where $\C=\add\C_0$, the smallest additive subcategory
containing $\C_0$ and closed under direct summands.
If $c\in\C$, $\sph n c$, $\UP(n,c)$ and so on mean
$\sph n \add c$, $\UP(n,\add c)$ and so on.
A $\C_0$-exact sequence means an $\add\C_0$-exact sequence.
A sequence $\Bbb E$ in $\A$ is $\C_0$-exact if and only if
for any $c\in\C_0$, $\A(\Bbb E,c)$ is exact.

\paragraph
By definition, any object of $\C$ is an injective object in $\A_{\C}$.

\paragraph
Let $\E$ be an exact category, and $\I$ an additive subcategory of $\E$.
Then for $e\in \E$, we define 
\begin{multline*}
\Push_\E(n,\I):=\{e\in \E\mid \text{There exists an exact sequence}\\
0\rightarrow e\rightarrow c^0\rightarrow
c^1\rightarrow\cdots\rightarrow
c^{n-1} \text{ with $c^i\in \I$}\}.
\end{multline*}
Note that $\Push_\E(0,\I)$ is the whole $\E$.
Thus $\Push_{\A_\C}(n,\C)=\UP_{\A}(n,\C)$.

If $a\in\E$ is a direct summand of an object of $\I$, then $a\in\Push(\infty,\I)$.

\begin{lemma}\label{UP-short.lem}
  Let $\E$ be an exact category.
  Let $\I$ be an additive subcategory of $\E$
  consisting of injective objects.
  Let 
  \[
  0\rightarrow a \xrightarrow f a' \xrightarrow g a'' \rightarrow 0
  \]
  be an exact sequence in $\E$ and $m\geq 0$.
  Then
  \begin{enumerate}
  \item[\bf 1] If $a\in\Push(m,\I)$ and $a''\in\Push(m,\I)$, then
    $a'\in\Push(m,\I)$.
  \item[\bf 2] If $a'\in \Push(m+1,\I)$ and $a''\in\Push(m,\I)$, then
    $a\in\Push(m+1,\I)$.
  \item[\bf 3] If $a\in\Push(m+1,\I)$, $a'\in\Push(m,\I)$, then
    $a''\in\Push(m,\I)$.
  \end{enumerate}
\end{lemma}

\begin{proof}
Let $i:\E\hookrightarrow\Cal A$ be the Gabriel--Quillen embedding \cite{TT}.
We consider that $\E$ is a full subcategory of $\Cal A$ closed under
extensions, and a sequence in $\E$ is exact if and only if it is so in $\Cal A$.

We prove {\bf 1}.
We use induction on $m$.
The case that $m=0$ is trivial, and so we assume that $m>0$.
Let
\[
0\rightarrow a\rightarrow c\rightarrow b\rightarrow 0
\]
be an exact sequence such that $c\in \I$ and $b\in\Push(m-1,\I)$.
Let
\[
0\rightarrow a''\rightarrow c''\rightarrow b''\rightarrow 0
\]
be an exact sequence such that $c''\in \I$ and $b''\in\Push(m-1,\I)$.
As $\C(a',c)\rightarrow \C(a,c)$ is surjective, we can form a
commutative diagram with exact rows and columns
\[
\xymatrix{
  & 0 \ar[d] & 0 \ar[d] & 0 \ar[d] \\
  0 \ar[r]
  & a \ar[r]^f \ar[d] & a' \ar[r]^g \ar[d] & a'' \ar[r] \ar[d] & 0 \\
  0 \ar[r] & c \ar[r]^{\scriptsize\begin{pmatrix}1\\0
  \end{pmatrix}} \ar[d] & c\oplus c'' \ar[r]^-{
   (1\;0)
  } \ar[d] & c'' \ar[r] \ar[d] & 0 \\
  0 \ar[r] & b \ar[r] \ar[d] & b' \ar[r] \ar[d] & b'' \ar[r] \ar[d] & 0 \\
  & 0 & 0 & 0
}
\]
in $\A$.
As $\E$ is closed under extensions in $\A$, this diagram is a diagram in $\E$.
By induction assumption, $b'\in \Push(m-1,\I)$.
Hence $a'\in \Push(m,\I)$.

We prove {\bf 2}.
Let $0\rightarrow a'\rightarrow c\rightarrow b'\rightarrow 0$ be an exact sequence
in $\E$ such that $c\in\I$ and $b'\in\Push(m,\I)$.
Then we have a commutative diagram in $\E$ with exact rows and columns
\[
\xymatrix{
  & & 0 \ar[d] & 0 \ar[d] \\
  & 0 \ar[r] & a \ar[r]^f \ar[d] & a' \ar[d] \ar[r]^g & a'' \ar[r] & 0 \\
  & 0 \ar[r] & c \ar[r]^{1_c} \ar[d] & c \ar[d] \ar[r] & 0 \\
  0 \ar[r] & a'' \ar[r] & b \ar[d] \ar[r] & b' \ar[d] \ar[r] & 0 \\
  & & 0 & 0
}.
\]
Applying {\bf 1}, which we have already proved, $b\in\Push(m,\I)$, since
$a''$ and $b'$ lie in $\Push(m,\I)$.
So $a\in\Push(m+1,\I)$, as desired.

We prove {\bf 3}.
Let $0\rightarrow a\rightarrow c\rightarrow b\rightarrow 0$ be an exact sequence
in $\E$ such that $c\in \I$ and $b\in\Push(m,\I)$.
Taking the push-out diagram
\[
\xymatrix{
& 0 \ar[d] & 0 \ar[d] \\
0 \ar[r] & a \ar[r]^f \ar[d] & a' \ar[r]^g \ar[d] & a'' \ar[r] \ar[d]^{1_{a''}} & 0 \\
0 \ar[r] & c \ar[d] \ar[r] & u \ar[d] \ar[r] & a'' \ar[r] & 0 \\
& b \ar[r]^{1_b}  \ar[d] & b \ar[d] \\
& 0 & 0
}.
\]
Then $u\in\Push(m,\I)$ by {\bf 1}, which we have already proved.
Since $c\in I$, the middle row splits.
Then by the exact sequence $0\rightarrow a''\rightarrow u\rightarrow c\rightarrow 0$ and {\bf 2},
we have that $a''\in\Push(m,\I)$, as desired.
\end{proof}

\begin{corollary}
  Let $\E$ and $\I$ be as in Lemma~\ref{UP-short.lem}.
  Let $m\geq 0$, and $a,a'\in\E$.
  Then $a\oplus a'\in\Push(m,\I)$ if and only if $a,a'\in\Push(m,\I)$.
\end{corollary}

\begin{proof}
  The \lq if' part is obvious by Lemma~\ref{UP-short.lem}, {\bf 1},
  considering the exact sequence
  \begin{equation}\label{split.eq}
  0\rightarrow a\rightarrow a\oplus a'\rightarrow a'\rightarrow 0.
  \end{equation}

  We prove the \lq only if' part by induction on $m$.
  If $m=0$, then there is nothing to prove.
  Let $m>0$.
  Then by induction assumption, $a'\in\Push(m-1,\I)$.
  Then applying Lemma~\ref{UP-short.lem}, {\bf 2} to the
  exact sequence (\ref{split.eq}), we have that $a\in\Push(m,\I)$.
  $a'\in\Push(m,\I)$ is proved similarly.
\end{proof}

\begin{corollary}\label{UP-short.cor}
  Let 
  \[
  0\rightarrow a \xrightarrow f a' \xrightarrow g a'' \rightarrow 0
  \]
  be a $\C$-exact sequence in $\A$ and $m\geq 0$.
  Then
  \begin{enumerate}
  \item[\bf 1] If $a\in\UP(m,\C)$ and $a''\in\UP(m,\C)$, then
    $a'\in\UP(m,\C)$.
  \item[\bf 2] If $a'\in \UP(m+1,\C)$ and $a''\in\UP(m,\C)$, then
    $a\in\UP(m+1,\C)$.
  \item[\bf 3] If $a\in\UP(m+1,\C)$, $a'\in\UP(m,\C)$, then
    $a''\in\UP(m,\C)$.
  \end{enumerate}
\qed
\end{corollary}

\paragraph
We define
$\sph{}\C=\sph \infty \C:=\bigcap_{i\geq 0}\sph i \C$ and
  $\UP(\infty,\C):=\bigcap_{j\geq 0}\UP(j,\C)$.
Obviously, $\C\subset \UP(\infty,\C)$.

\begin{lemma}
  We have
  \begin{multline*}
  \UP(\infty,\C)=\{a\in\A\mid \text{There exists some $\C$-exact
    sequence }\\
  0\rightarrow a \rightarrow c^0\rightarrow c^1\rightarrow
  c^2\rightarrow \cdots \text{ with $c^i\in\C$ for $i\geq 0$}\}.
  \end{multline*}
\end{lemma}

\begin{proof}
  Let $a\in\UP(\infty,\C)$, and
  take any $\C$-exact sequence
  \[
  0\rightarrow a\rightarrow c^0\rightarrow a^1\rightarrow 0
  \]
  with $c^0\in\C$.
  Then $a^1\in \UP(\infty,\C)$ by
  Corollary~\ref{UP-short.cor}, and we can continue infinitely.
\end{proof}

\paragraph
We define $\Y_{\infty}:=\bigcup_{i\geq 0} \Y_i$.
We also define $\X_{i,j}:=\sph i \C\cap \UP(j,\C)$ for $0\leq i,j\leq\infty$.

\paragraph
Let $0\leq i,j\leq\infty$.
We say that $a\in\A$ lies in $\Z_{i,j}$ if
there is a short exact sequence
\[
0\rightarrow y\rightarrow x \rightarrow a \rightarrow 0
\]
in $\A$ such that $x\in \X_{i,j}$ and $y\in \Y_i$.

\paragraph We define $\infty\pm r = \infty$ for $r\in\Bbb R$.

\begin{lemma}\label{descending-Taka.lem}
  Let $0\leq i,j\leq\infty$ with $j\geq 1$.
  Assume that $\C\subset \sph{i+1}\C$.
  Let $0\rightarrow z \xrightarrow f x \xrightarrow g z'\rightarrow 0$ be
  a short exact sequence in $\A$ with $z\in\Z_{i,j}$ and
  $x\in\X_{i+1,j-1}$.
  Then $z'\in\Z_{i+1,j-1}$.
\end{lemma}

\begin{proof}
  By assumption, there is an exact sequence
  \[
  0\rightarrow y \xrightarrow j x'\xrightarrow \varphi z\rightarrow 0
  \]
  such that $\CCdim y < i$ and $x'\in \X_{i,j}$.
  As $j\geq 1$, there is an $\C$-exact sequence
  \[
  0\rightarrow x'\xrightarrow{h} c \rightarrow x'''\rightarrow 0
  \]
  such that $c\in \C$.
  Then we have a commutative diagram with exact rows and columns
  \[
  \xymatrix{
    & 0 \ar[d] & 0 \ar[d] & 0 \ar[d] \\
    0 \ar[r] & y \ar[r]^{hj} \ar[d]^j
    & c \ar[d]^{\scriptsize\begin{pmatrix} 1 \\ 0 
    \end{pmatrix}} \ar[r] & y' \ar[r] \ar[d] & 0 \\
    0 \ar[r] & x' \ar[d]^{\varphi} \ar[r]^{\scriptsize\begin{pmatrix} h \\ f\varphi 
      \end{pmatrix}} & c \oplus x \ar[d]^{\scriptsize\begin{pmatrix} 0 \;1 
    \end{pmatrix}} \ar[r] & x'' \ar[d] \ar[r] & 0 \\
    0 \ar[r] & z \ar[r]^f \ar[d] & x \ar[r]^g \ar[d] & z' \ar[r] \ar[d]
    & 0 \\
    & 0 & 0 & 0
  }
  \]
  As the top row is exact, $y\in \Y_i$, and $c\in\C$, $y'\in \Y_{i+1}$.
  By assumption, $c\in \X_{i+1,\infty}$ and $x\in\X_{i+1,j-1}$.
  So $c\oplus x\in\X_{i+1,j-1}$.
  As the middle row is $\C$-exact and $x'\in \X_{i,j}$,
  we have that $x''\in\X_{i+1,j-1}$ by Corollary~\ref{UP-short.cor}.
  The right column shows that $z'\in\Z_{i+1,j-1}$, as desired.
\end{proof}

\begin{lemma}\label{ascending-Taka.lem}
  Let $0\leq i,j\leq \infty$, and assume that $i\geq 1$ and
  $\C\subset \sph{i}\C$.
  Let
  \begin{equation}\label{Takahashi-short.eq}
    0\rightarrow z\xrightarrow f x \xrightarrow g z'\rightarrow 0
  \end{equation}
    be
  a short exact sequence in $\A$ with $z'\in \Z_{i,j}$ and $x\in\X_{i,j+1}$.
  Then $z\in\Z_{i-1,j+1}$.
\end{lemma}

\begin{proof}
  Take an exact sequence $0\rightarrow y'\rightarrow x''\xrightarrow h z'
  \rightarrow 0$ such that $x''\in\X_{i,j}$ and $y'\in\Y_i$.
  Taking the pull-back of (\ref{Takahashi-short.eq}) by $h$,
  we get a commutative diagram with exact rows and columns
  \[
  \xymatrix{
    & & 0 \ar[d] & 0 \ar[d] \\
    & 0 \ar[r] \ar[d] & y' \ar[r]^{1_{y'}} \ar[d] &
    y' \ar[r] \ar[d] & 0\\
    0 \ar[r] & z \ar[d]^{1_z} \ar[r]^j & a \ar[r] \ar[d] & x'' \ar[r]
    \ar[d]^h & 0 \\
    0 \ar[r] & z \ar[r]^f \ar[d] & x \ar[r]^g \ar[d] &
    z' \ar[r] \ar[d] & 0 \\
    & 0 & 0 & 0
  }.
  \]
  By induction, we can prove easily that $\sph i \C\subset \sph{i+1-l}\Y_l$.
  In particular, $\sph i\C\subset \sph 1 \Y_i$,
  and $\Ext_{\A}^1(x,y')=0$.
  Hence the middle column splits, and we can replace $a$ by $x\oplus y'$.
  By the definition of $\Y_i$, there is an exact sequence
  \[
  0\rightarrow y\rightarrow c \rightarrow y'\rightarrow 0
  \]
  of $\A$ such that $y\in \Y_{i-1}$ and $c\in\C$.
  Then adding $1_x$ to this sequence, we get
  \[
  0\rightarrow y\rightarrow x\oplus c\rightarrow x\oplus y'\rightarrow 0
  \]
  is exact.
  Pulling back this exact sequence with $j:z\rightarrow a= x\oplus y'$,
  we get a commutative diagram with exact rows and columns
  \[
  \xymatrix{
    & 0 \ar[d] & 0 \ar[d] & 0 \ar[d] \\
    0 \ar[r] & y \ar[r] \ar[d]^{1_y} & x' \ar[d]\ar[r]& z \ar[r] \ar[d]^j
    & 0 \\
    0 \ar[r] & y \ar[d]\ar[r] & x\oplus c \ar[d] \ar[r] &
    x\oplus y'\ar[r]\ar[d] & 0 \\
    & 0 \ar[r] & x'' \ar[r]^{1_{x''}} \ar[d] & x'' \ar[r] \ar[d] & 0
    \\
    & & 0 & 0
  }.
  \]
  As $x''\in \sph 1\C$, the middle column is $\C$-exact.
  As $x''\in\X_{i,j}$ and $x\oplus c\in \X_{i,j+1}$, we have that
  $x'\in \X_{i-1,j+1}$.
  As the top row shows, $z\in \Z_{i-1,j+1}$, as desired.
\end{proof}

\begin{theorem}\label{X_{m,n}-approximation.thm}
  Let $0\leq n,m\leq \infty$, and assume that $\C\subset\sph n \C$.
  For $z\in \A$, the following are equivalent.
  \begin{enumerate}
  \item[\bf 1] $z\in \Z_{n,m}$.
  \item[\bf 2] There is an exact sequence
    \begin{equation}\label{x-resol.eq}
      0\rightarrow x_n\xrightarrow{d_n}
      x_{n-1}\xrightarrow{d_{n-1}} x_0\xrightarrow \varepsilon z\rightarrow
    0
    \end{equation}
    such that $x_i\in \X_{n-i,m+i}$.
  \end{enumerate}
  If, moreover, for each $a\in \A$, there is a surjection
  $x\rightarrow a$ with $x\in \X_{n,n+m}$, then these conditions are
  equivalent to the following.
  \begin{enumerate}
  \item[\bf 3] For each exact sequence {\rm(\ref{x-resol.eq})} with
    $x_i\in \X_{n-i,m+i+1}$ for $0\leq i\leq n-1$, we have that
    $x_n\in\X_{0,n+m}$.
  \end{enumerate}
\end{theorem}

\begin{proof}
  {\bf 1$\Rightarrow$2}.
  There is an exact sequence $0\rightarrow y\rightarrow x_0\xrightarrow
  \varepsilon z   \rightarrow 0$ with $x_0\in \X_{n,m}$ and $y\in \Y_n$.
  So there is an exact sequence
  \[
  0\rightarrow x_n\xrightarrow{d_n} x_{n-1}\xrightarrow{d_{n-1}}
  \cdots\xrightarrow{d_2} x_1\rightarrow y\rightarrow 0
  \]
  with $x_i\in \C$ for $1\leq i\leq n$.
  As $\C\subset \X_{n,\infty}$, we are done.

  {\bf 2$\Rightarrow$1}.
  Let $z_i=\Image d_i$ for $i=1,\ldots,n$, and $z_0:=z$.
  Then by descending induction on $i$,
  we can prove $z_i\in \Z_{n-i,m+i}$ for $i=n,n-1,\ldots,0$, using
  Lemma~\ref{descending-Taka.lem} easily.

  {\bf 1$\Rightarrow$3} is also proved easily, using
  Lemma~\ref{ascending-Taka.lem}.

  {\bf 3$\Rightarrow$2} is trivial.
\end{proof}
  
\section{$(n,C)$-TF property}\label{(n,C)-TF.sec}
\label{takahashi.sec}

\paragraph
In the rest of this paper, let
$\Lambda$ be a module-finite $R$-algebra, which may not be commutative.
A $\Lambda$-bimodule means a $\Lambda\otimes_R\Lambda\op$-module.
Let $C\in \md \Lambda$ be fixed.
Set $\Gamma:=\End_{\Lambda\op} C$.
Note that $\Gamma$ is also a module-finite $R$-algebra.
We denote $(?)^\dagger:=\Hom\Lop(?,C):\md\Lambda\rightarrow (\Gamma\md)\op$,
and $(?)^\ddagger:=\Hom_\Gamma(?,C):\Gamma\md\rightarrow (\md\Lambda)\op$.

\paragraph
We denote
$\Syz_{\md \Lambda}(n,C)$,
$\UP_{\md \Lambda}(n,C)$,
and
$\Cdim_{\md\Lambda} M$
respectively
by $\Syz\Lop(n,C)$, $\UP\Lop(n,C)$, and $\Cdim\Lop M$.

\paragraph
Note that for $M\in\md\Lambda$ and $N\in\Gamma\md$,
we have standard isomorphisms
\begin{equation}\label{standard-isom.eq}
  \Hom\Lop(M,N^\ddagger)\cong \Hom_{\Gamma\otimes_R \Lambda\op}(N\otimes_R M,C)
  \cong \Hom_\Gamma(N,M^\dagger).
\end{equation}
The first isomorphism sends $f:M\rightarrow N^\ddagger$ to
the map $(n\otimes m\mapsto f(m)(n))$.
Its inverse is given by $g:N\otimes_R M\rightarrow C$ to
$(m\mapsto (n\mapsto g(n\otimes
m)))$.
This shows that $(?)^\dagger$ has $((?)^\ddagger)\op:(\Gamma\md)\op
\rightarrow\md \Lambda$ as a right adjoint.
Hence $((?)^\dagger)\op$ is right adjoint to $(?)^\ddagger$.
We denote the unit of adjunction
$\Id\rightarrow (?)^{\dagger\ddagger}=(?)^\ddagger(?)^\dagger$
by $\lambda$.
Note that for $M\in \md\Lambda$, the map $\lambda_M:
M\rightarrow M^{\dagger\ddagger}$ is given by
$\lambda_M(m)(\psi)=\psi(m)$ for $m\in M$ and
$\psi\in M^\dagger=\Hom\Lop(M,C)$.
We denote the unit of adjunction $N\rightarrow N^{\ddagger\dagger}$ by
$\mu=\mu_N$ for
$N\in\Gamma\md$.
When we view $\mu$ as a morphism $N^{\ddagger\dagger}\rightarrow N$ (in the opposite
category $(\Gamma\md)\op$), then it is the counit of adjunction.

\begin{lemma}
  $(?)^\dagger$ and $(?)^\ddagger$ give a contravariant
  equivalence between
  $\add C\subset \md \Lambda$ and $\add\Gamma\subset \Gamma\md$.
\end{lemma}

\begin{proof}
  It suffices to show that $\lambda:M\rightarrow M^{\dagger\ddagger}$ is an isomorphism for
  $M\in\add C$, and $\mu:N\rightarrow N^{\ddagger\dagger}$ is an isomorphism for
  $N\in\add\Gamma$.
  To verify this, we may assume that $M=C$ and $N=\Gamma$.
  This case is trivial.
\end{proof}

\begin{definition}[cf.~{\cite[(2.2)]{Takahashi}}]\label{(n,C)-TF.def}
  Let $M\in\md\Lambda$.
  We say that $M$ is $(1,C)$-TF or $M\in\TF\Lop(1,C)$ if $\lambda_M
  :M\rightarrow M^{\dagger
    \ddagger}$ is injective.
  We say that $M$ is $(2,C)$-TF or $M\in\TF\Lop(2,C)$ if $\lambda_M:
  M\rightarrow M^{\dagger
    \ddagger}$ is bijective.
  Let $n\geq 3$.
  We say that $M$ is $(n,C)$-TF or $M\in \TF\Lop(n,C)$
  if $M$ is $(2,C)$-TF and
  $\Ext^i_\Gamma(M^\dagger,C)=0$ for $1\leq i\leq n-2$.
  As a convention, we define that any $M\in\md \Lambda$ is $(0,C)$-TF.
\end{definition}

\begin{lemma}\label{Takahashi-prelimilary.lem}
  Let $\Theta: 0\rightarrow M\rightarrow L\rightarrow N\rightarrow 0$ be
  a $C$-exact sequence in $\md \Lambda$.
  Then for $n\geq 0$, we have the following.
  \begin{enumerate}
  \item[\bf 1] If $M\in\TF(n,C)$ and $N\in\TF(n,C)$, then $L\in\TF(n,C)$.
  \item[\bf 2] If $L\in \TF(n+1,C)$ and $N\in\TF(n,C)$, then $M\in\TF(n+1,C)$.
  \item[\bf 3] If $M\in \TF(n+1,C)$ and $L\in\TF(n,C)$, then $N\in\TF(n,C)$.
  \end{enumerate}
\end{lemma}

\begin{proof}
  We have a commutative diagram
\[
\xymatrix{
  0\ar[r] & M \ar[d]^{\lambda_M} \ar[r]^h & L \ar[d]^{\lambda_L} \ar[r] &
  N \ar[d]^{\lambda_N} \ar[r] & 0 \\
  0\ar[r] & M\ddd \ar[r]^{h\ddd} & L\ddd \ar[r] &
  N\ddd \ar[r] & \Ext^1_\Gamma(M\da,C)\ar[r] & \Ext^1_\Gamma(L\da,C) \ar[r] & \cdots
}
\]
with exact rows.

We only prove {\bf 3}.
We may assume that $n\geq 1$.
So $\lambda_M$ is an isomorphism and $\lambda_L$ is injective.
By the five lemma, 
$\lambda_N$ is injective, and the case that $n=1$ has been done.
If $n\geq 2$, then $\lambda_L$ is also an isomorphism and
$\Ext^1_\Gamma(M\da,C)=0$, and so $\lambda_N$ is an isomorphism.
Moreover, for $1\leq i\leq n-2$, 
$\Ext^i_\Gamma(L\da,C)$ and $\Ext^{i+1}_\Gamma(M\da,C)$ vanish.
so $\Ext^i_\Gamma(N\da,C)=0$ for $1\leq i\leq n-2$, and hence $N\in\TF(n,C)$.

{\bf 1} and {\bf 2} are also proved similarly.
\end{proof}

\begin{lemma}[cf.~{\cite[Proposition~3.2]{Takahashi}}]\label{Takahashi.lem}
  \begin{enumerate}
  \item[\bf 1] For $n=0,1$, $\Syz\Lop(n,C)=\UP\Lop(n,C)$.
  \item[\bf 2]     For $n\geq 0$, $\TF\Lop(n,C)=\UP\Lop(n,C)$.
  \end{enumerate}
\end{lemma}
  
\begin{proof}
  If $n=0$, then $\Syz\Lop(n,C)
  =\TF\Lop(0,C)=\UP\Lop(0,C)=\md \Lambda$.
  So we may assume that $n\geq 1$.

  Let $M\in\Syz\Lop(1,C)$.
  Then there is an injection $\varphi:M\rightarrow N$ with $N\in\add C$.
Then
\[
\xymatrix{
  M \ar[d]^{\lambda_M} \ar@{^(->}[r]^\varphi & N \ar[d]^{\lambda_{N}}_\cong \\
  M^{\dagger\ddagger} \ar[r]^{\varphi\ddd} & N^{\dagger\ddagger}
}
\]
is a commutative diagram.
So $\lambda_M$ is injective, and $M\in\TF\Lop(1,C)$.
This shows $\UP\Lop(1,C)\subset\Syz\Lop(1,C)\subset \TF\Lop(1,C)$.
So {\bf 2$\Rightarrow$1}.

We prove {\bf 2}.
First, we prove $\UP\Lop(n,C)\subset \TF\Lop(n,C)$ for $n\geq 1$.
We use induction on $n$.
The case $n=1$ is already done above.

Let $n\geq 2$ and $M\in\UP\Lop(n,C)$.
Then by the definition of $\UP\Lop(n,C)$,
there is a $C$-exact sequence
\[
0\rightarrow M\rightarrow L\rightarrow N\rightarrow 0
\]
such that $L\in\add C$ and $N\in\UP\Lop(n-1,C)$.
By induction hypothesis, $N\in\TF\Lop(n-1,C)$.
Hence $M\in\TF\Lop(n,C)$ by Lemma~\ref{Takahashi-prelimilary.lem}.
We have proved that $\UP\Lop(n,C)\subset \TF\Lop(n,C)$.

Next we show that $\TF\Lop(n,C)\subset \UP\Lop(n,C)$ for $n\geq 1$.
We use induction on $n$.

Let $n=1$.
Let $\rho: F\rightarrow M\da$ be any surjective $\Gamma$-linear map with
$F\in\add\Gamma$.
Then the map $\rho':M\rightarrow F\dda$ which corresponds to $\rho$ by
the adjunction (\ref{standard-isom.eq}) is
\[
\rho':M\xrightarrow{\lambda_M}M\ddd \xrightarrow{\rho\dda}F\dda,
\]
which is injective by assumption.
Then $\rho$ is the composite
\[
\rho: F \xrightarrow{\mu_F} F\dddd \xrightarrow{(\rho')\da}M\da,
\]
which is a surjective map by assumption.
So $(\rho')\da$ is also surjective, and hence
$\rho':M\rightarrow F\dda$ gives a $(1,C)$-universal pushforward.

Now let $n\geq 2$.
By what we have proved, $M$ has a $(1,C)$-universal pushforward
$h:M\rightarrow L$.
Let $N=\Coker h$.
Then we have a $C$-exact sequence
\[
0\rightarrow M\rightarrow L\rightarrow N\rightarrow 0
\]
with $L\in\add C$.
As $M\in \TF(n,C)$, $N\in\TF(n-1,C)$
by Lemma~\ref{Takahashi-prelimilary.lem}.
By induction assumption, $N\in\UP(n-1,C)$.
So by the definition of $\UP(n,C)$, we have that $M\in \UP(n,C)$, as desired.
\end{proof}

\begin{lemma}\label{Syz(2,C).lem}
  For any $N\in\Gamma\md$, $N\dda\in\Syz(2,C)$.
\end{lemma}

\begin{proof}
  Let
  \[
  F_1\xrightarrow h F_0\rightarrow N \rightarrow 0
  \]
  be an exact sequence in $\Gamma\md$ such that $F_i\in\add\Gamma$.
  Then
  \[
  0\rightarrow N\dda\rightarrow F_0\dda\xrightarrow{h\dda}F_1\dda
  \]
  is exact, and $F_i\dda\in \add C$.
  This shows that $N\dda\in\Syz(2,C)$.
\end{proof}

\paragraph
We denote by $(S_n')_C=(S_n')_C^{\Lambda\op,R}$
the class of $M\in \md\Lambda$ such that
$M$ viewed as an $R$-module lies in $(S_n')_C^R$, see (\ref{(S_n')_N.par}).

\begin{lemma}\label{Syz-S_r'.lem}
  Assume that $C$ satisfies $(S_n')$ as an $R$-module.
  Then $\Syz(r,C)\subset (S_r')_C^{\Lambda\op,R}$ for $r\geq 1$.
\end{lemma}

\begin{proof}
  This follows easily from Corollary~\ref{S_n-increase.cor}.
\end{proof}

\paragraph
For an additive category $\C$ and its additive subcategory $\X$, we denote
by $\C/\X$ the quotient of $\C$ divided by the ideal consisting of
morphisms which factor through objects of $\X$.

\paragraph
For each $M\in\md \Lambda$, take a presentation
\begin{equation}\label{M-presentation.eq}
  \Bbb F(M): F_1(M)\xrightarrow{\partial}F_0(M)\xrightarrow{\varepsilon}
  M\rightarrow 0
\end{equation}
with $F_i\in\add \Lambda_\Lambda$.
We denote
\[
\Coker(\partial\da)=\Coker(1_C\otimes \partial^t)=C\otimes_\Lambda \Tr M
\]
by $\Tr_C M$, where $(?)^t=\Hom\Lop(?,\Lambda)$ and $\Tr$ is the transpose,
see \cite[(V.2)]{ASS}, and we call it the $C$-transpose of $M$.
$\Tr_C$ is an additive
functor from $\smd \Lambda:=\md\Lambda/\add\Lambda_\Lambda$
to $\Gamma {}_C\smd:=\Gamma\md/\add C$.

\begin{proposition}
  Let $n\geq 0$, and assume that $\Ext^i_\Gamma(C,C)=0$ for $i=1,\ldots,n$.
  Then for $M\in\md \Lambda$, we have the following.
  \begin{enumerate}
  \item[\bf 0] For $1\leq i\leq n$, $\Ext^i_\Gamma(\Tr_C ?,C)$
    is a well-defined additive functor $\smd\Lambda\rightarrow \md\Lambda$.
  \item[\bf 1] 
    If $n=1$, there is an exact sequence
    \[
    0\rightarrow \Ext^1_\Gamma(\Tr_CM,C)
    \rightarrow M\xrightarrow{\lambda_M}M\ddd
    \rightarrow \Ext^2_\Gamma(\Tr_C M,C).
    \]
    If $n=0$, then there is an injective homomorphism $\Ker\lambda_M\hookrightarrow
    \Ext^1_\Gamma(\Tr_CM,C)$.
  \item[\bf 2] If $n\geq 2$, then
    \begin{enumerate}
    \item[\bf i] There is an exact sequence
      \[
    0\rightarrow \Ext^1_\Gamma(\Tr_CM,C)
    \rightarrow M\xrightarrow{\lambda_M}M\ddd
    \rightarrow \Ext^2_\Gamma(\Tr_C M,C)\rightarrow 0.
    \]
  \item[\bf ii] There are isomorphisms
    $\Ext^{i+2}\Gamma(\Tr_CM,C) \cong \Ext^i_\Gamma(M\da,C)$ for
    $1\leq i\leq n-2$.
  \item[\bf iii] There is an injective map
    $\Ext^{n-1}_\Gamma(M\da,C)\hookrightarrow \Ext^{n+1}_\Gamma(\Tr_C M,C)$.
    \end{enumerate}
  \end{enumerate}
\end{proposition}

\begin{proof}
  {\bf 0} is obvious by assumption.

  We consider that $\Bbb F(M)$ is a complex with $M$ at degree zero.
  Then consider
  \[
  \Bbb Q(M):=\Bbb F(M)\da[2]:
  F_1(M)\da\xleftarrow{\partial\da} F_0(M)\da \xleftarrow{\varepsilon\da}
  M\da \leftarrow 0
  \]
  where $F_1(M)\da$ is at degree zero.
  As this complex is quasi-isomorphic to $\Tr_C(M)$, there is a spectral
  sequence
  \[
  E_1^{p,q}=\Ext^q_\Gamma(\Bbb Q(M)^{-p},C)\Rightarrow \Ext_\Gamma^{p+q}(\Tr_CM,C).
  \]
  In general, $\Ker\lambda_M=E_2^{1,0}\cong E_\infty^{1,0}\subset E^1$.
  If $n\geq 1$, then $E_1^{0,1}=0$, and $E_\infty^{1,0}=E^1$.
  Moreover, as $E_1^{0,1}=0$, 
  $\Coker\lambda_M\cong E_2^{2,0}\cong E_\infty^{2,0}\subset E^2$.
  So {\bf 1} follows.

  If $n\geq 2$, then $E_1^{0,2}=E_1^{1,1}=0$ by assumption, so
  $E_\infty^{2,0}=E^2$, and {\bf i} of {\bf 2} follows.
  Note that $E_1^{p,q}=0$ for
  $p\geq 3$.
  Moreover, $E_1^{p,q}=0$ for $p=0,1$ and $1\leq q\leq n$.
  So for $1\leq i\leq n-1$, we have
  \[
  E_1^{2,i}\cong E_\infty^{2,i}\hookrightarrow E^{i+2},
  \]
  and the inclusion is an isomorphism if $1\leq i\leq n-2$.
  So {\bf ii} and {\bf iii} of {\bf 2} follow.
\end{proof}
  
\begin{corollary}
  Let $n\geq 1$.
  If $\Ext^i_\Gamma(C,C)=0$ for $1\leq i\leq n$, then $M$ is $(n,C)$-TF
  if and only if $\Ext^i_\Gamma(\Tr_C M,C)=0$ for $1\leq i\leq n$.
  If $\Ext^i_\Gamma(C,C)=0$ for $1\leq i < n$ and
  $\Ext^i_\Gamma(\Tr_CM,C)=0$ for $1\leq i\leq n$, then $M$ is $(n,C)$-TF.
\end{corollary}

\section{Canonical module}\label{canonical.sec}

\paragraph
Let $R=(R,\fm)$ be semilocal, where $\fm$ is the
Jacobson radical of $R$.

\paragraph
We say that a dualizing complex $\Bbb I$ over $R$ is {\em normalized} if for any
maximal ideal $\fn$ of $R$, $\Ext^0_R(R/\fn,\Bbb I)\neq 0$.
We follow the definition of \cite{Hartshorne}.

\paragraph
For a left or right $\Lambda$-module $M$, 
$\dim M$ or $\dim_\Lambda M$
denotes the dimension $\dim_R M$ of $M$,
which is independent of the choice of $R$.
We call $\depth_R(\fm,M)$, which is also independent of $R$,
the global depth, $\Lambda$-depth, or depth of $M$,
and denote it by $\depth_\Lambda M$ or $\depth M$.
$M$ is called globally Cohen--Macaulay or GCM for short, 
if $\dim M=\depth M$.
$M$ is GCM
if and only if it is Cohen--Macaulay as an $R$-module,
and all the maximal ideals of $R$ have the same height.
This notion is independent of $R$, and depends only on
$\Lambda$ and $M$.
$M$ is called a globally
maximal Cohen--Macaulay (GMCM for short) if
$\dim \Lambda=\depth M$.
We say that the algebra $\Lambda$ is GCM
if the $\Lambda$-module $\Lambda$ is GCM.
However, in what follows, if $R$ happens to be local, then
GCM and Cohen--Macaulay 
(resp.\ GMSM and maximal Cohen--Macaulay)
(over $R$)
are the same thing, and used interchangeably.

\paragraph\label{canonical-complete-def.par}
Assume that $(R,\fm)$ is complete semilocal, and $\Lambda\neq 0$.
Let $\Bbb I$ be a normalized dualizing complex of $R$.
The lowest non-vanishing cohomology group $\Ext^{-s}_R(\Lambda,\Bbb I)$
($\Ext^i_R(\Lambda,\Bbb I)=0$ for $i<-s$) is denoted by $K_\Lambda$, and is called the
{\em canonical module} of $\Lambda$.
Note that $K_\Lambda$ is a $\Lambda$-bimodule.
Hence it is also a $\Lambda\op$-bimodule.
In this sense, $K_\Lambda=K_{\Lambda\op}$.
If $\Lambda=0$, then we define $K_\Lambda=0$.

\paragraph
Let $S$ be the center of $\Lambda$.
Then $S$ is module-finite over $R$, and $\Bbb I_S=\RHom_R(S,\Bbb I)$ is a
normalized dualizing complex of $S$.
This shows that $\RHom_R(\Lambda,\Bbb I)\cong \RHom_S(\Lambda,\Bbb I_S)$, and
hence the definition of $K_\Lambda$ is also independent of $R$.

\begin{lemma}
  The number $s$ in {\rm(\ref{canonical-complete-def.par})}
  is nothing but $d:=\dim \Lambda$.
  Moreover, 
\[
\Ass_R K_\Lambda=\Assh_R \Lambda:=\{P\in\Min_R \Lambda\mid \dim R/P=\dim \Lambda\}.
\]
\end{lemma}

\begin{proof}
  We may replace $R$ by $R/\ann_R \Lambda$, and may assume that $\Lambda$ is a
  faithful module.
  We may assume that $\Bbb I$ is a fundamental dualizing complex of $R$.
  That is, for each $P\in\Spec R$, $E(R/P)$, the injective hull of $R/P$, appears
  exactly once (at dimension $-\dim R/P$).
  If $\Ext^{-i}_R(\Lambda,\Bbb I)\neq 0$, then there exists some $P\in\Spec R$ such that
  $\Ext^{-i}_{R_P}(\Lambda_P,\Bbb I_P)\neq 0$.
  Then $P\in\Supp_R\Lambda$ and $\dim R/P\geq i$.
  On the other hand, $\Ext^{-d}_{R_P}(\Lambda_P,\Bbb I_P)$ has length $l(\Lambda_P)$
  and is nonzero for $P\in \Assh_R \Lambda$.
  So $s=d$.

  The argument above shows that each $P\in\Assh_R\Lambda=\Assh R$ supports $K_\Lambda$.
  So $\Assh_R\Lambda\subset\Min_R K_\Lambda$.
  On the other hand, as the complex $\Bbb I$ starts at degree $-d$,
  $K_\Lambda\subset \Bbb I^{-d}$, and $\Ass K_\Lambda\subset \Ass\Bbb I^{-d}\subset
  \Assh R=\Assh_R\Lambda$.
\end{proof}

  \begin{lemma}
Let $(R,\fm)$ be complete semilocal.
Then $K_\Lambda$ satisfies the $(S_2^\Lambda)^R$-condition.
  \end{lemma}

  \begin{proof}
    It is easy to see that $(K_\Lambda)_\fn$ is either zero or $K_{\Lambda_\fn}$ for
    each maximal ideal $\fn$ of $R$.
    Hence we may assume that $R$ is local.
    Replacing $R$ by $R/\ann_R \Lambda$, we may assume that $\Lambda$ is a
    faithful $R$-module, and we are to prove that $K_\Lambda$ satisfies
    $(S_2')^R$ by Lemma~\ref{S_n-S_n'.lem}.
    Replacing $R$ by a Noether normalization, we may further assume that $R$ 
    is regular by Lemma~\ref{finite-s_n.lem}, {\bf 1}.
    Then $K_\Lambda=\Hom_R(\Lambda,R)$.
    So $K_\Lambda\in\Syz(2,R)\subset (S_2')^R$
    by Lemma~\ref{Syz(2,C).lem}
    (consider that $\Lambda$ there is $R$ here,
    and $C$ there is also $R$ here).
  \end{proof}

  \paragraph
Assume that $(R,\fm)$ is semilocal which may not be complete.
We say that a finitely generated $\Lambda$-bimodule $K$ is
a {\em canonical module} of $\Lambda$
if $\hat K$ is isomorphic to
the canonical module $K_{\hat \Lambda}$ as a 
$\hat \Lambda$-bimodule.
It is unique up to isomorphisms, and denoted by $K_\Lambda$.
We say that $K \in \md \Lambda$ is a right canonical module of $\Lambda$
if $\hat K$ is isomorphic to $K_{\hat \Lambda}$ in $\md \hat \Lambda$,
where $\hat ?$ is the $\fm$-adic completion.
If $K_\Lambda$ exists, then $K$ is a right canonical module if and only if
$K\cong K_\Lambda$ in $\md\Lambda$.

These definitions are independent of $R$, in the sense that the (right) canonical
module over $R$ and that over the center of $\Lambda$ are the same thing.
The right canonical module of $\Lambda\op$ is called the left canonical module.
A $\Lambda$-bimodule $\omega$ is said to be a weakly canonical bimodule if
${}_\Lambda\omega$ is left canonical, and $\omega_\Lambda$ is right canonical.
The canonical module $K_{\Lambda\op}$ 
of $\Lambda\op$ is canonically identified with $K_\Lambda$.

\paragraph\label{canonical-localization.par}
If $R$ has a normalized
dualizing complex $\Bbb I$, then $\hat{\Bbb I}$ is a normalized dualizing complex
of $\hat R$, and so it is easy to see that $K_\Lambda$ exists and
agrees with $\Ext^{-d}(\Lambda,\Bbb I)$, where $d=\dim \Lambda(:=\dim_R\Lambda)$.
In this case, for any $P\in\Spec R$, $\Bbb I_P$ is a dualizing complex of $R_P$.
So if $R$ has a dualizing complex and $(K_\Lambda)_P\neq 0$, then $(K_\Lambda)_P$,
which is the lowest nonzero cohomology group of $\RHom_{R_P}(\Lambda_P,\Bbb I_P)$, is
the $R_P$-canonical module of $\Lambda_P$.
See also Theorem~\ref{Aoyama.thm} below.

\begin{lemma}\label{non-complete.lem}
  Let $(R,\fm)$ be local, and assume that $K_\Lambda$ exists.
  Then we have the following.
  \begin{enumerate}
  \item[\bf 1] $\Ass_R K_\Lambda=\Assh_R \Lambda$.
  \item[\bf 2] $K_\Lambda\in (S_2^\Lambda)^{R}$.
  \item[\bf 3] $R/\ann K_\Lambda$ is quasi-unmixed, and hence is universally catenary.
  \end{enumerate}
\end{lemma}

\begin{proof}
  All the assertions are proved easily using the case that $R$ is complete.
\end{proof}

\paragraph
A $\Lambda$-module $M$ is said to be
$\Lambda$-full over $R$ if $\Supp_R M =\Supp_R \Lambda$.

\begin{lemma}\label{ogoma.lem}
  Let $(R,\fm)$ be local.
  If $K_\Lambda$ exists and $\Lambda$ satisfies the $(S_2)^R$-condition, then
  $R/\ann K_\Lambda$ is equidimensional, and $K_\Lambda$ is $\Lambda$-full
  over $R$.
\end{lemma}

\begin{proof}
  The same as
  the proof of \cite[Lemma~4.1]{Ogoma} (use Lemma~\ref{non-complete.lem}, {\bf 3}).
\end{proof}

\paragraph\label{local-duality.par}
Let $(R,\fm)$ be local, and $\Bbb I$ be a normalized dualizing complex.
By the local duality,
\[
K_\Lambda^\vee=\Ext^{-d}(\Lambda,\Bbb I)^\vee\cong H^d_\fm(\Lambda)
\]
(as $\Lambda$-bimodules),
where $E_R(R/\fm)$ is the injective hull of the $R$-module $R/\fm$,
and $(?)^\vee$ is the Matlis dual $\Hom_R(?,E_R(R/\fm))$.

\paragraph\label{contravariant-equivalence.par}
Let $(R,\fm)$ be semilocal, and $\Bbb I$ be a normalized dualizing complex.
Note that $\RHom_R(?,\Bbb I)$ induces a contravariant equivalence between
$D\fg(\Lambda\op)$ and $D\fg(\Lambda)$.
Let $\Bbb J\in D\fg(\Lambda\otimes_R \Lambda\op)$ be $\RHom_R(\Lambda,\Bbb I)$.
\[
\RHom_R(?,\Bbb I):D\fg(\Lambda\op)\rightarrow D\fg(\Lambda)
\]
is identified with
\[
\RHom\Lop(?,\RHom_R({}_\Lambda\Lambda_R,\Bbb I))=\RHom\Lop(?,\Bbb J)
\]
and similarly,
\[
\RHom_R(?,\Bbb I):D\fg(\Lambda)\rightarrow D\fg(\Lambda\op)
\]
is identified with
$\RHom_\Lambda(?,\Bbb J)$.
Note that a left or right $\Lambda$-module $M$
is maximal Cohen--Macaulay if and only if $\RHom_R(M,\Bbb I)$ is concentrated
in degree $-d$, where $d=\dim \Lambda$.

\paragraph $\Bbb J$ above is a dualizing complex of $\Lambda$ in the sense
of Yekutieli \cite[(3.3)]{Yekutieli}.

\paragraph\label{ext-vanishing.par}
$\Lambda$ is GCM if and only if $K_\Lambda[d]\rightarrow
\Bbb J$ is an isomorphism.
If so, $M\in \md \Lambda$ is GMCM if and only if
$\RHom_R(M,\Bbb I)$ is concentrated in degree $-d$ if and only if
$\Ext\Lop^i(M,K_\Lambda)=0$ for $i>0$.
Also, in this case,
as $K_\Lambda[d]$ is a dualizing complex, it is of finite injective dimension
both as a left and a right $\Lambda$-module.
To prove these, we may take the completion, and may assume that
$R$ is complete.
All the assertions are independent of $R$, so taking the Noether normalization,
we may assume that $R$ is local.
By (\ref{contravariant-equivalence.par}), the assertions follow.

\paragraph\label{ext-dual.par}
For any $M\in\md \Lambda$ which is GMCM,
\[
M
\cong
\RHom_R(\RHom_R(M,\Bbb I),\Bbb I)
\cong
\RHom_R(\Ext^{-d}\Lop(M,K_\Lambda[d]),\Bbb I)[-d].
\]
Hence $M\da:=\Hom\Lop(M,K_\Lambda)$
is also a GMCM $\Lambda$-module, 
and hence
\[
\Hom_{\Lambda}(M\da,K_\Lambda)\rightarrow \RHom_{\Lambda}(M\da,\Bbb J)=\RHom_R(M\da,\Bbb I)
\]
is an isomorphism (in other words, $\Ext_\Lambda^i(M\da,K_\Lambda)=0$ for $i>0$).
So the canonical map 
\begin{equation}\label{double-canonical.eq}
M\rightarrow \Hom_\Lambda(\Hom\Lop(M,K_\Lambda),K_\Lambda)=\Hom_\Lambda(M\da,K_\Lambda)
\end{equation}
$m\mapsto (\varphi\mapsto \varphi m)$
is an isomorphism.
This isomorphism is true without assuming that $R$ has a dualizing complex
(but assuming the existence of a canonical module), passing to the completion.
Note that if $\Lambda=R$ and $K_R$ exists and Cohen--Macaulay, then
$K_R$ is a dualizing complex of $R$.

Similarly, for $N\in\Lambda\md$ which is GMCM,
\[
N\rightarrow \Hom\Lop(\Hom_\Lambda(N,K_\Lambda),K_\Lambda)
\]
$n\mapsto (\varphi\mapsto \varphi n)$
is an isomorphism.

\paragraph\label{algebra-isom.par}
In particular, letting $M=\Lambda$, if $\Lambda$ is GCM,
we have that $K_\Lambda=\Hom\Lop(\Lambda,K_\Lambda)$ is GMCM.
Moreover, 
\[
\Lambda\rightarrow \End\Lop K_\Lambda
\]
is an $R$-algebra isomorphism, where $a\in\Lambda$ goes to the left multiplication by
$a$.
Similarly,
\[
\Lambda\rightarrow (\End_\Lambda K_\Lambda)\op
\]
is an isomorphism of $R$-algebras.

\paragraph\label{canonical-top-local-cohomology.par}
Let $(R,\fm)$ be a $d$-dimensional complete local ring, and $\dim \Lambda=d$.
Then by the local duality,
\[
H^d_\fm(K_\Lambda)^\vee
\cong
\Ext^{-d}_R(K_\Lambda,\Bbb I)
\cong
\Ext^{-d}\Lop(K_\Lambda,\Bbb J)
\cong \End\Lop K_\Lambda,
\]
where $\Bbb J=\Hom_R(\Lambda,\Bbb I)$ and $(?)^\vee=\Hom_R(?,E_R(R/\fm))$.

\section{$n$-canonical module}\label{n-canonical.sec}

\paragraph
We say that $\omega$
is an $R$-semicanonical right $\Lambda$-module (resp.\
$R$-semicanonical left $\Lambda$-module,
weakly $R$-semicanonical $\Lambda$-bimodule,
$R$-semicanonical $\Lambda$-bimodule) if for any $P\in\Spec R$,
$R_P \otimes_R \omega$ is
the right canonical module (resp.\ left canonical module, 
weakly canonical module, canonical module)
of $R_P\otimes_R \Lambda$ 
for any $P\in\supp_R \omega$.
If we do not mention what $R$ is, then it may mean $R$ is the center of $\Lambda$.
An $R$-semicanonical right $\Lambda\op$-module (resp.\
$R$-semicanonical left $\Lambda\op$-module,
weakly $R$-semicanonical $\Lambda\op$-bimodule,
$R$-semicanonical $\Lambda\op$-bimodule) is nothing but
an $R$-semicanonical left $\Lambda$-module (resp.\
$R$-semicanonical right $\Lambda$-module,
weakly $R$-semicanonical $\Lambda$-bimodule,
$R$-semicanonical $\Lambda$-bimodule).

\paragraph
Let $C\in\md \Lambda$ (resp.\ $\Lambda\md$, $(\Lambda\otimes_R \Lambda\op)\md$,
$(\Lambda\otimes_R \Lambda\op)\md$).
We say that $C$ is an
$n$-canonical right $\Lambda$-module (resp.\
$n$-canonical left $\Lambda$-module,
weakly $n$-canonical $\Lambda$-bimodule,
$n$-canonical $\Lambda$-bimodule)
over $R$ if $C\in (S_n')^R$,
and for each $P\in R\an{<n}$, we have that $C_P$ is an
$R_P$-semicanonical right $\Lambda_P$-module (resp.\
$R_P$-semicanonical left $\Lambda_P$-module,
weakly $R_P$-semicanonical $\Lambda_P$-bimodule,
$R_P$-semicanonical $\Lambda_P$-bimodule).
If we do not mention what $R$ is, it may mean $R$ is the center of $\Lambda$.

\begin{example}\label{n-canonical.ex}
  \begin{enumerate}
    \item[\bf 0] The zero module $0$ is an $R$-semicanonical $\Lambda$-bimodule.
    \item[\bf 1]
  If $R$ has a dualizing complex $\Bbb I$, then the lowest non-vanishing
  cohomology group $K:=\Ext^{-s}_R(\Lambda,\Bbb I)$ is an $R$-semicanonical
  $\Lambda$-bimodule.
\item[\bf 2] 
  By Lemma~\ref{non-complete.lem}, any left or right
  $R$-semicanonical module $K$ of $\Lambda$
  satisfies the $(S_2^\Lambda)^R$-condition.
  Thus a (right) semicanonical module is
  $2$-canonical over $R/\ann_R\Lambda$.
\item[\bf 3] If $K$ is (right) semicanonical (resp.\ $n$-canonical)
  and $L$ is a
  projective $R$-module such that $L_P$ is rank at most one,
  then $K\otimes_R L$ is again (right) semicanonical (resp.\ $n$-canonical).
\item[\bf 4] If $R$ is a normal domain and $C$ its rank-one reflexive module of $R$,
  then $C$ is a $2$-canonical $R$-module (here $\Lambda=R$).
\item[\bf 5]  The $R$-module $R$ is $n$-canonical if and only if
  for $P\in R\bra{<n}$, $R_P$ is Gorenstein.
  This is equivalent to say that $R$ satisfies
  $(T_{n-1})+(S_n)$.
  \end{enumerate}
\end{example}

\paragraph
As in section~\ref{takahashi.sec},
let $C\in\md\Lambda$, and set
$\Gamma=\End\Lop C$, $(?)\da=\Hom\Lop(?,C)$, and $(?)\dda=\Hom_\Gamma(?,C)$.
Moreover, we set
$\Lambda_1:=(\End_\Gamma C)\op$.
The $R$-algebra map $\Psi_1:\Lambda\rightarrow\Lambda_1$ is induced by the
right action of $\Lambda$ on $C$.

\begin{lemma}\label{1-canonical.lem}
  Let $C\in\md\Lambda$ be a $1$-canonical $\Lambda\op$-module over $R$.
  Let $M\in\md \Lambda$.
  Then the following are equivalent.
  \begin{enumerate}
  \item[\bf 1] $M\in \TF(1,C)$.
  \item[\bf 2] $M\in \UP(1,C)$.
  \item[\bf 3] $M\in\Syz(1,C)$.
  \item[\bf 4] $M\in (S_1')^R_C$.
  \end{enumerate}
\end{lemma}

\begin{proof}
  {\bf 1$\Leftrightarrow$2} is Lemma~\ref{Takahashi.lem}.
  {\bf 2$\Rightarrow$3} is trivial.
  {\bf 3$\Rightarrow$4} follows from Lemma~\ref{Syz-S_r'.lem} immediately.

  We prove {\bf 4$\Rightarrow$1}.
  We want to prove that $\lambda_M: M\rightarrow M\ddd$ is injective.
  By Example~\ref{injective-bijective.ex}, localizing at each $P\in R\an0$,
  we may assume that $(R,\fm)$ is zero-dimensional local.
  We may assume that $M$ is nonzero.
  By assumption, $C$ is nonzero, and hence $C=K_\Lambda$ by assumption.
  As $R$ is zero-dimensional, $\Lambda$ is GCM, and hence
  $\Lambda\rightarrow \Gamma=\End\Lop K_\Lambda$ is an isomorphism by
  (\ref{algebra-isom.par}).
  As $\Lambda$ is GCM and
  $M$ is GMCM, (\ref{double-canonical.eq}) is
  an isomorphism.
  As $\Lambda=\Gamma$, the result follows.
\end{proof}

\begin{lemma}\label{happy.lem}
  Let $C$ be a $1$-canonical right $\Lambda$-module over $R$,
  and $N\in\Gamma\md$.
  Then $N\dda\in \TF\Lop(2,C)$.
  Similarly, for $M\in\md\Lambda$, we have that $M\da\in\TF_{\Gamma}(2,C)$.
\end{lemma}

\begin{proof}
  Note that $\lambda_{N\dda}:N\dda\rightarrow N^{\ddagger\dagger\ddagger}$ is a
  split monomorphism.
  Indeed, $(\mu_N)\dda:N^{\ddagger\dagger\ddagger}\rightarrow N\dda$ is the left inverse.
  Assume that
  $N\dda\notin\TF(2,C)$, then $W:=\Coker \lambda_{N\dda}$ is nonzero.
  Let $P\in\Ass_R W$.
  As $W$ is a submodule of $N^{\ddagger\dagger\ddagger}$,
  $P\in\Ass_R N^{\ddagger\dagger\ddagger}\subset \Ass_R C\subset \Min R$.
  So $C_P$ is the right canonical module $K_{\Lambda_P}$.
  So $\Gamma_P=\Lambda_P$, and $(\lambda_{N\dda})_P$ is an isomorphism.
  This shows that $W_P=0$, and this is a contradiction.
  The second assertion is proved similarly.
\end{proof}

\begin{lemma}\label{Lambda-Lambda_1-injective.lem}
  Let $(R,\fm)$ be local, and assume that $K_\Lambda$ exists.
  Let $C:=K_\Lambda$.
  If $\Lambda$ is GCM, $\Psi_1:\Lambda\rightarrow \Lambda_1$
  is an isomorphism.
\end{lemma}

\begin{proof}
  As $C$ possesses a bimodule structure, we have a canonical map
  $\Lambda\rightarrow\Gamma=\End\Lop C$,
  which is an isomorphism as $\Lambda$ is GCM by (\ref{algebra-isom.par}).
  So $\Lambda_1$ is identified with $\Delta=(\End_\Lambda C)\op$.
  Then $\Psi_1:\Lambda\rightarrow(\End_\Lambda C)\op$ is an isomorphism
  again by (\ref{algebra-isom.par}).
\end{proof}

\begin{lemma}
  If $C$ satisfies the $(S_1')^R$ condition, then $\Gamma\in (S_1')^R_C$
  and $\Lambda_1\in (S_1')^R_C$.
  Moreover, $\Ass_R \Gamma =\Ass_R \Lambda_1=\Ass_R C=\Min_R C$.
\end{lemma}

\begin{proof}
  The first assertion is by $\Gamma=\Hom\Lop(C,C)\in\Syz_\Gamma(2,C)$,
  and $\Lambda_1=\Hom_\Gamma(C,C)=\Syz_{\Lambda_1}(2,C)$.
  We prove the second assertion.
  $\Ass_R\Gamma\subset\Ass_R \End_R C=\Ass_R C$.
  $\Ass_R\Lambda_1\subset \Ass_R \End_R C=\Ass_R C=\Min_R C$.
  It remains to show that
  $\Supp_R C=\Supp_R\Gamma=\Supp_R \Lambda_1$.
  Let $P\in\Spec R$.
  If $C_P=0$, then $\Gamma_P=0$ and $(\Lambda_1)_P=0$.
  On the other hand, if $C_P\neq 0$, then the identity map
  $C_P\rightarrow C_P$ is not zero, and hence
  $\Gamma_P\neq 0$ and $(\Lambda_1)_P\neq 0$.
\end{proof}

\paragraph
Let $C$ be a $1$-canonical right $\Lambda$-module over $R$.
Define $Q:=\prod_{P\in\Min_R C}R_P$.
If $P\in\Min_R C$, then $C_P=K_{\Lambda_P}$.
Hence $\Phi_P:\Lambda_P\rightarrow(\Lambda_1)_P$ is an
isomorphism by Lemma~\ref{Lambda-Lambda_1-injective.lem}.
So $1_Q\otimes \Psi_1:
Q\otimes_R \Lambda\rightarrow Q\otimes_R \Lambda_1$ is also an isomorphism.
As $\Ass_R\Lambda_1=\Min_R C$, we have that $\Lambda_1\subset Q\otimes_R \Lambda_1$.

\begin{lemma}
Let $C$ be a $1$-canonical right $\Lambda$-module over $R$.
  If $\Lambda$ is commutative, then so are $\Lambda_1$ and $\Gamma$.
\end{lemma}

\begin{proof}
  As $\Lambda_1\subset Q\otimes_R\Lambda_1=Q\otimes_R\Lambda$ and
  $Q\otimes_R \Lambda$ is commutative, $\Lambda_1$ is a commutative ring.
  We prove that $\Gamma$ is commutative.
  As $\Ass_R \Gamma\subset\Min_R C$, $\Gamma$ is a subring of $Q\otimes\Gamma$.
  As
  \[
  Q\otimes_R\Gamma\cong\prod_{P\in\Min_R C}\End_{\Lambda_P}C_P\cong
\prod_P \End_{\Lambda_P}(K_{\Lambda_P})
\]
and $\Lambda_P\rightarrow\End_{\Lambda_P}(K_{\Lambda_P})$ is an isomorphism
(as $\Lambda_P$ is zero-dimensional), $Q\otimes_R \Gamma$ is, and hence
$\Gamma$ is also, commutative.
\end{proof}

\begin{lemma}\label{Lambda-Lambda_1-full.lem}
  Let $C$ be a $1$-canonical right $\Lambda$-module over $R$.
  Let $M$ and $N$ be left \(resp.\ right, bi-\) modules of $\Lambda_1$,
  and assume that $N\in (S_1')^{\Lambda_1,R}$.
  Let $\varphi:M\rightarrow N$ be a $\Lambda$-homomorphism of
  left \(resp.\ right, bi-\) modules.
  Then $\varphi$ is a $\Lambda_1$-homomorphism of left \(resp.\ right, bi-\)
  modules.
\end{lemma}

\begin{proof}
  Let $Q=\prod_{P\in\Min_R C}R_P$.
  Then we have a commutative diagram
  \[
  \xymatrix{
    M \ar[r]^\varphi \ar[d]^{i_M} & N \ar[d]^{i_N} \\
    Q\otimes_R M \ar[r]^{1\otimes\varphi} & Q\otimes_R N
  },
  \]
  where $i_M(m)=1\otimes m$ and $i_N(n)=1\otimes n$.
  Clearly, $i_M$ and $i_N$ are $\Lambda_1$-linear.
  As $\varphi$ is $\Lambda$-linear,
  $1\otimes\varphi$ is $Q\otimes\Lambda$-linear.
  Since $\Lambda_1\subset Q\otimes\Lambda_1=Q\otimes\Lambda$,
  $1\otimes\varphi$ is $\Lambda_1$-linear.
  As $i_N$ is injective, it is easy to see that
  $\varphi$ is $\Lambda_1$-linear.
\end{proof}

\begin{lemma}\label{fully-faithful.lem}
  Let $C$ be a $1$-canonical right $\Lambda$-module over $R$.
  Then the restriction $M\mapsto M$ is a full and faithful functor
  from $(S_1')^{\Lambda_1,R}$ to $(S_1')_C^{\Lambda,R}$.
  Similarly, it gives a full and faithful functors
  $(S_1')^{\Lambda_1\op,R}\rightarrow (S_1')_C^{\Lambda\op,R}$ and
  $(S_1')^{\Lambda_1\otimes_R \Lambda_1\op,R}\rightarrow (S_1')_C^{\Lambda\otimes_R
    \Lambda\op,R}$.
\end{lemma}

\begin{proof}
  We only consider the case of left modules.
  If $M\in \Lambda_1\md$, then it is a homomorphic image of
  $\Lambda_1 \otimes_R M$.
  Hence $\supp_R M\subset \supp_R \Lambda_1\subset\supp_R C$.
  So the functor is well-defined and obviously faithful.
  By Lemma~\ref{Lambda-Lambda_1-full.lem}, it is also full, and we are done.
\end{proof}

\paragraph
Let $C$ be a $1$-canonical $\Lambda$-bimodule over $R$.
Then the left action of $\Lambda$ on $C$ induces an $R$-algebra map
$\Phi:\Lambda\rightarrow \Gamma=\End\Lop C$.
Let $Q=\prod_{P\in\Min_R C}R_P$.
Then $\Gamma\subset Q\otimes_R \Gamma = Q\otimes_R \Lambda$.
From this we get

\begin{lemma}\label{Lambda-Gamma-full.lem}
  Let $C$ be a $1$-canonical $\Lambda$-bimodule over $R$.
  Let $M$ and $N$ be left \(resp.\ right, bi-\) modules of $\Gamma$,
  and assume that $N\in (S_1')^{\Gamma,R}$.
  Let $\varphi:M\rightarrow N$ be a $\Lambda$-homomorphism of
  left \(resp.\ right, bi-\) modules.
  Then $\varphi$ is a $\Gamma$-homomorphism of left \(resp.\ right, bi-\)
  modules.
\end{lemma}

\begin{proof}
  Similar to Lemma~\ref{Lambda-Lambda_1-full.lem}, and left to the reader.
\end{proof}

\begin{corollary}\label{ddd-dastar.cor}
  Let $C$ be as above.
  $(?)\ddd=\Hom_\Gamma(\Hom\Lop(?,C),C)$ is canonically isomorphic to
  $(?)^{\dagger\star}=\Hom_\Lambda(\Hom\Lop(?,C),C)$, where $(?)^\star=\Hom_\Lambda(?,C)$.
\end{corollary}

\begin{proof}
This is immediate by Lemma~\ref{Lambda-Gamma-full.lem}.
\end{proof}

\begin{lemma}\label{Lambda-Gamma-fully-faithful.lem}
Let $C$ be a $1$-canonical $\Lambda$-bimodule over $R$.
Then $\Phi$ induces a full and faithful functor
$(S_1')^{\Gamma,R}\rightarrow (S_1')^{\Lambda,R}_C$.
Similarly, 
$(S_1')^{\Gamma\op,R}\rightarrow (S_1')^{\Lambda\op,R}_C$
and
$(S_1')^{\Gamma\otimes_R\Gamma\op,R}\rightarrow (S_1')^{\Lambda\otimes_R\Lambda\op,R}_C$
are also induced.
\end{lemma}

\begin{proof}
  Similar to Lemma~\ref{fully-faithful.lem}, and left to the reader.
\end{proof}

\begin{corollary}\label{Lambda_1=Delta.cor}
  Let $C$ be a $1$-canonical $\Lambda$-bimodule.
  Set $\Delta:=(\End_\Lambda C)\op$.
  Then the canonical map $\Lambda\rightarrow\Gamma$ induces
  an equality
  \[
  \Lambda_1=(\End_\Gamma C)\op = (\End_{\Lambda}C)\op = \Delta.
\]
Similarly, we have
\[
\Lambda_2:=\End_{\Delta\op}C=\End\Lop C=\Gamma.
\]
\end{corollary}

\begin{proof}
  As $C\in (S_1')^{\Gamma,R}$, the first assertion follows from
  Lemma~\ref{Lambda-Gamma-fully-faithful.lem}.
  The second assertion is proved by left-right symmetry.
\end{proof}

\begin{lemma}\label{1-canonical-injective-right.lem}
  Let $C$ be a $1$-canonical right $\Lambda$-module over $R$.
  Set $\Lambda_1:=(\End_\Gamma C)\op$.
  Let $\Psi_1:\Lambda\rightarrow \Lambda_1$ be the canonical map induced by the
  right action of $\Lambda$ on $C$.
  Then $\Psi_1$ is injective if and only if $\Lambda$ satisfies the
  $(S_1')^R$ condition and $C$ is $\Lambda$-full over $R$.
\end{lemma}

\begin{proof}
  $\Psi_1:\Lambda\rightarrow \Lambda_1$ is nothing but $\lambda_\Lambda:
  \Lambda\rightarrow\Lambda\ddd$, and the result follows from
  Lemma~\ref{1-canonical.lem} immediately.
\end{proof}

\begin{lemma}\label{1-canonical-injective.lem}
  Let $C$ be a $1$-canonical $\Lambda$-bimodule over $R$.
  Then the following are equivalent.
  \begin{enumerate}
  \item[\bf 1] The canonical map $\Psi:\Lambda\rightarrow\Delta$ is injective,
    where $\Delta=(\End_\Lambda C)\op$, and the map is induced by the right
    action of $\Lambda$ on $C$.
  \item[\bf 2] $\Lambda$ satisfies the $(S_1')^R$ condition, and $C$ is
    $\Lambda$-full over $R$.
  \item[\bf 3] The canonical map $\Phi: \Lambda\rightarrow\Gamma$ is
    injective, where the map is induced by the left action of $\Lambda$ on $C$.
  \end{enumerate}
\end{lemma}

\begin{proof}
    By Corollary~\ref{Lambda_1=Delta.cor}, 
    we have that $\Lambda_1 =(\End_\Gamma C)\op =\Delta$.
    So {\bf 1$\Leftrightarrow$2} is a consequence of
    Lemma~\ref{1-canonical-injective-right.lem}.

    Reversing the roles of the left and the right, we get
    {\bf 2$\Leftrightarrow$\bf 3} immediately.
\end{proof}

\begin{lemma}\label{C-bimod-isom.lem}
  Let $C$ be a $1$-canonical right $\Lambda$-module over $R$.
  Then the canonical map
  \begin{equation}\label{lambda-gamma.eq}
  \Hom\Lop(\Lambda_1,C)\rightarrow \Hom\Lop(\Lambda,C)\cong C
  \end{equation}
  induced by the canonical map $\Psi_1:\Lambda\rightarrow \Lambda_1$ is an
  isomorphism of $\Gamma\otimes_R\Lambda_1\op$-modules.
\end{lemma}

\begin{proof}
  The composite map
  \[
  C\cong\Hom_{\Lambda_1}(\Lambda_1,C)
  =
  \Hom_\Lambda(\Lambda_1,C)
  \rightarrow
  \Hom_{\Lambda}(\Lambda,C)
  \cong C
  \]
  is the identity.
  The map is a $\Gamma\otimes_R\Lambda\op$-homomorphism.
  It is also $\Lambda_1\op$-linear by Lemma~\ref{fully-faithful.lem}.
\end{proof}

\paragraph When $(R,\fm)$ is local and $C=K_\Lambda$, then $\Lambda_1=\Delta$,
and the map
(\ref{lambda-gamma.eq}) is an isomorphism of $\Gamma\otimes_R\Delta\op$-modules
from $K_\Delta$ and $K_\Lambda$,
where $\Delta=(\End_\Lambda K_\Lambda)\op$.
Indeed, to verify this, we may assume that $R$ is complete regular local with
$\ann_R \Lambda=0$, and hence $C=\Hom_R(\Lambda,R)$, and
$C$ is a $2$-canonical $\Lambda$-bimodule over $R$, see (\ref{n-canonical.ex}).
So (\ref{Lambda_1=Delta.cor}) and Lemma~\ref{C-bimod-isom.lem} apply.
Hence we have

\begin{corollary}\label{lambda-gamma2.cor}
  Let $(R,\fm)$ be a local ring with a canonical module $C=K_\Lambda$ of $\Lambda$.
  Then $K_\Delta=\Hom\Lop(\Delta,K_\Lambda)$
  is isomorphic to $K_\Lambda$ as a $\Gamma\otimes_R\Delta\op$-module,
  where $\Delta=(\End_\Lambda K_\Lambda)\op$.
  \qed
\end{corollary}

\begin{lemma}\label{n-canonical-left-right.lem}
  Let $n\geq 1$.
  If $C$ is an $n$-canonical right $\Lambda$-module over $R$, then
  \begin{enumerate}
  \item[\bf 1] $C$ is an $n$-canonical right $\Lambda_1$-module over $R$.
  \item[\bf 2] $C$ is an $n$-canonical left $\Gamma$-module over $R$.
  \end{enumerate}
\end{lemma}

\begin{proof}
  {\bf 1}.
  As the $(S_n')$-condition holds, it suffices to prove that for $P\in R\an{<n}$,
  $C_P\cong (K_{\Lambda_1})_P$ as a right $(\Lambda_1)_P$-module.
  After localization, replacing $R$ by $R_P$, we may assume that
  $R$ is local and $C=K_\Lambda$.
  Then $C\cong K_{\Lambda}\cong K_{\Lambda_1}$ as right $\Lambda$-modules.
  Both $C$ and $K_{\Lambda_1}$ are in $(S_1')^{\Lambda_1\op,R}$, and isomorphic
  in $\md\Lambda$.
  So they are isomorphic in $\md\Lambda_1$ by Lemma~\ref{fully-faithful.lem}.

  {\bf 2}.
  Similarly, assuming that $R$ is local and $C=K_\Lambda$, it suffices to show that
  $C\cong K_\Gamma$ as left $\Gamma$-modules.
  Identifying $\Gamma=\End_{\Delta\op}C=\Lambda_2$ and using the left-right symmetry,
  this is the same as the proof of {\bf 1}.
\end{proof}

\begin{lemma}\label{2-canonical.lem}
  Let $C\in\md\Lambda$ be a $2$-canonical right $\Lambda$-module over $R$.
  Let $M\in\md \Lambda$.
  Then the following are equivalent.
  \begin{enumerate}
  \item[\bf 1] $M\in \TF(2,C)$.
  \item[\bf 2] $M\in \UP(2,C)$.
  \item[\bf 3] $M\in\Syz(2,C)$.
  \item[\bf 4] $M\in (S_2')_C^R$.
  \end{enumerate}
\end{lemma}

\begin{proof}
  We may assume that $\Lambda$ is a faithful $R$-module.
  {\bf 1$\Leftrightarrow$2$\Rightarrow$3$\Rightarrow$4} is easy.
  We show {\bf 4$\Rightarrow$1}.
  By Example~\ref{injective-bijective.ex}, localizing at each $P\in R\an{\leq 1}$,
  we may assume that $R$ is a Noetherian local ring of dimension at most one.
  So the formal fibers of $R$ are zero-dimensional, and hence
  $\hat M\in (S_2')_{\hat C}^{\hat R}$, where $\hat ?$ denotes the completion.
  So we may further assume that $R=(R,\fm)$ is complete local.
  We may assume that $M\neq 0$ so that $C\neq 0$ and hence $C=K_\Lambda$.
  The case $\dim R=0$ is similar to the proof of Lemma~\ref{1-canonical.lem},
  so we prove the case that $\dim R=1$.
  Note that $I=H^0_\fm(\Lambda)$ is a two-sided ideal of $\Lambda$,
  and any module in $(S_1')^{\Lambda\op,R}$ is annihilated by $I$.
  Replacing $\Lambda$ by $\Lambda/I$,
  we may assume that $\Lambda$ is a 
  maximal Cohen--Macaulay $R$-module.
  Then (\ref{double-canonical.eq}) is an isomorphism.
  As $C=K_\Lambda$ and
  \[
  \Lambda\rightarrow \End\Lop K_\Lambda=\End\Lop C=\Gamma
  \]
  is an $R$-algebra isomorphism, we have that $\lambda_M:M\rightarrow M\ddd$ is
  identified with the isomorphism (\ref{double-canonical.eq}), as desired.
\end{proof}

\begin{corollary}\label{2-canonical-isom.cor}
  Let $C$ be a $2$-canonical right $\Lambda$-module over $R$.
  Then the canonical map
  $\Phi:\Lambda\rightarrow\Lambda_1$ is an isomorphism if and only if
  $\Lambda$ satisfies $(S_2')^R$ and $C$ is full.
\end{corollary}

\begin{proof}
  Follows immediately by Lemma~\ref{2-canonical.lem} applied to $M=\Lambda$.
\end{proof}

\paragraph\label{Gamma-Delta.par}
Let $C$ be a $2$-canonical $\Lambda$-bimodule.
Let $\Gamma=\End\Lop C$ and $\Delta=(\End_\Lambda C)\op$.
Then by the left multiplication, an $R$-algebra map $\Lambda\rightarrow \Gamma$
is induced, while by the right multiplication, an $R$-algebra map
$\Lambda\rightarrow \Delta$ is induced.
Let $Q=\prod_{P\in\Min_R C}R_P$.
Then as $\Gamma\subset Q\otimes_R \Gamma=Q\otimes_R\Lambda=Q\otimes_R\Delta\supset\Delta$,
both $\Gamma$ and $\Delta$ are identified with $Q$-subalgebras of $Q\otimes_R\Lambda$.
As $\Delta=\Lambda_1=\Lambda\ddd$, we have a commutative diagram
\[
\xymatrix{
  \Lambda \ar[r]^{\lambda_\Lambda} \ar[d]^\nu & \Lambda\ddd \ar[d]^{\nu\ddd} & = \Delta  \\
  \Gamma \ar[r]^{\lambda_\Gamma} & \Gamma\ddd
}.
\]
As $\Gamma=\Hom_{\Lop}(C,C)=C^\dagger$, $\Gamma\in\Syz_\Lambda(2,C)$ by
Lemma~\ref{Syz(2,C).lem}.
By Lemma~\ref{2-canonical.lem}, we have that $\Gamma\in (S_2')_C$.
Hence by Lemma~\ref{2-canonical.lem} again, 
$\lambda_\Gamma:\Gamma\rightarrow \Gamma\ddd$ is an isomorphism.
Hence $\Delta\subset \Gamma$.
By symmetry $\Delta\supset\Gamma$.
So $\Delta=\Gamma$.
With this identification, $\Gamma$ acts on $C$ not only from left, but also
from right.
As the actions of $\Gamma$ extend those of $\Lambda$, $C$ is a $\Gamma$-bimodule.
Indeed, for $a\in\Lambda$, the left multiplication $\lambda_a:C\rightarrow C$
($\lambda_a(c)=ac$) is right $\Gamma$-linear.
So for $b\in\Gamma$, $\rho_b:C\rightarrow C$ ($\rho_b(c)=cb$) is
left $\Lambda$-linear,
and hence is left $\Gamma$-linear.

\begin{theorem}\label{2-canonical-equiv.thm}
  Let $C$ be a $2$-canonical right $\Lambda$-module.
  Then the restriction $M\mapsto M$ gives an equivalence
  $\rho: (S_2')^{\Lambda_1\op,R}\rightarrow (S_2')_C^{\Lambda\op,R}$.
\end{theorem}

\begin{proof}
  The functor is obviously well-defined, and is full and faithful by
  Lemma~\ref{fully-faithful.lem}.
  On the other hand, given $M\in(S_2')_C^{\Lambda\op,R}$, we have that
  $\lambda_M : M\rightarrow M\ddd$ is an isomorphism.
  As $M\ddd$ has a $\Lambda_1\op$-module
  structure which extends the $\Lambda\op$-module
  structure of $M\cong M\ddd$, we have that $\rho$ is also dense,
  and hence is an equivalence.
\end{proof}

\begin{corollary}\label{2-canonical-equiv.cor}
  Let $C$ be a $2$-canonical $\Lambda$-bimodule.
  Then the restriction $M\mapsto M$ gives an equivalence
  \[
  \rho: (S_2')_C^{\Gamma\otimes_R \Gamma\op,R}\rightarrow (S_2')_C^{\Lambda\otimes_R
    \Lambda\op,R}.
  \]
\end{corollary}

\begin{proof}
  $\rho$ is well-defined, and is obviously faithful.
  If $h:M\rightarrow N$ is a morphism of $(S_2)_C^{\Lambda\otimes_R \Lambda\op,R}$
  between objects of $(S_2)_C^{\Gamma\otimes_R\Gamma\op,R}$, then
  $h$ is $\Gamma$-linear $\Gamma\op$-linear by Theorem~\ref{2-canonical-equiv.thm}
  (note that $\Lambda_1=\Delta=\Gamma$ here).
  Hence $\rho$ is full.

  Let $M\in(S_2)_C^{\Lambda\otimes_R \Lambda\op,R}$, the left (resp.\ right)
  $\Lambda$-module
  structure of $M$ is extendable to that of a left (resp.\ right) $\Gamma$-module
  structure by Theorem~\ref{2-canonical-equiv.thm}.
  It remains to show that these structures make $M$ a $\Gamma$-bimodule.
  Let $a\in\Lambda$.
  Then $\lambda_a:M\rightarrow M$ given by $\lambda_a(m)=am$ is a right
  $\Lambda$-linear, and hence is right $\Gamma$-linear.
  So for $b\in\Gamma$, $\rho_b:M\rightarrow M$ given by $\rho_b(m)=mb$ is
  left $\Lambda$-linear, and hence is left $\Gamma$-linear, as desired.
\end{proof}

\begin{proposition}\label{2-canonical-duality.prop}
  Let $C$ be a $2$-canonical right $\Lambda$-module.
  Then $(?)\da: (S_2')_C^{\Lambda\op,R}\rightarrow (S_2')^{\Gamma,R}$ and
  $(?)\dda: (S_2')^{\Gamma,R}\rightarrow (S_2')_C^{\Lambda\op,R}$ give a
  contravariant equivalence.
\end{proposition}

\begin{proof}
  As we know that $(?)\da$ and $(?)\dda$ are contravariant adjoint each
  other, it suffices to show that the unit $\lambda_M:M\rightarrow M\ddd$
  and the (co-)unit $\mu_N: N\rightarrow N^{\ddagger\dagger}$ are isomorphisms.
  $\lambda_M$ is an isomorphism by Lemma~\ref{2-canonical.lem}.
  Note that $C$ is a $2$-canonical left $\Gamma$-module by
  Lemma~\ref{n-canonical-left-right.lem}.
  So $\mu_N$ is an isomorphism by Lemma~\ref{2-canonical.lem} applied to
  the right $\Gamma\op$-module $C$.
\end{proof}

\begin{corollary}
  Let $C$ be a $2$-canonical $\Lambda$-bimodule.
  Then $(?)\da=\Hom\Lop(?,C)$ and $\Hom_{\Lambda}(?,C)$ give a contravariant
  equivalence between $(S_2')_C^{\Lambda\op,R}$ and $(S_2')_C^{\Lambda,R}$.
  They also give a duality of $(S_2')_C^{\Lambda\otimes\Lambda\op,R}$.
\end{corollary}

\begin{proof}
  The first assertion
  is immediate by Proposition~\ref{2-canonical-duality.prop} and
  Theorem~\ref{2-canonical-equiv.thm}.
  The second assertion follows easily from the first and
  Corollary~\ref{2-canonical-equiv.cor}.
\end{proof}

\section{Non-commutative Aoyama's theorem}\label{Aoyama.sec}

\begin{lemma}\label{flat-local-free.lem}
  Let $(R,\fm,k)\rightarrow (R',\fm',k')$ be a flat local homomorphism
  between Noetherian local rings.
  \begin{enumerate}
  \item[\bf 1] Let $M$ be a $\Lambda$-bimodule such that $M':=R'\otimes_R M$
    is isomorphic to $\Lambda':=R'\otimes_R \Lambda$ as a $\Lambda'$-bimodule.
    Then $M\cong\Lambda$ as a $\Lambda$-bimodule.
  \item[\bf 2] Let $M$ be a right $\Lambda$ module such that $M':=R'\otimes_R M$
    is isomorphic to $\Lambda':=R'\otimes_R \Lambda$ as a right $\Lambda'$-module.
    Then $M\cong\Lambda$ as a right $\Lambda$-module.
  \end{enumerate}
\end{lemma}

\begin{proof}
  Taking the completion, we may assume that both $R$ and $R'$ are complete.
  Let $1=e_1+\cdots+e_r$ be the decomposition of $1$ into the mutually
  orthogonal primitive idempotents of the center $S$ of $\Lambda$.
  Then replacing $R$ by $Se_i$, $\Lambda$ by $\Lambda e_i$, and $R'$ by
  the local ring of $R'\otimes_R Se_i$ at any maximal ideal,
  we may further assume that $S=R$.
  This is equivalent to say that $R\rightarrow \End_{\Lambda\otimes_R \Lambda\op}\Lambda$
  is isomorphic.
  So $R'\rightarrow \End_{\Lambda'\otimes_{R'}(\Lambda')\op}\Lambda'$ is also isomorphic,
and hence the center of $\Lambda'$ is $R'$.

{\bf 1}.
Let $\psi:M'\rightarrow \Lambda'$ be an isomorphism.
Then we can write $\psi=\sum_{i=1}^m u_i\psi_i$ with $u_i\in R'$ and
$\psi_i\in \Hom_{\Lambda\otimes_R \Lambda\op}(M,\Lambda)$.
Also, we can write $\psi_i^{-1}=\sum_{j=1}^n v_j\varphi_j$ with $v_j\in R'$ and
$\varphi_j\in \Hom_{\Lambda\otimes_R \Lambda\op}(\Lambda,M)$.
As $\sum_{i,j}u_iv_j \psi_i\varphi_j=\psi \psi^{-1}=1\in\End_{\Lambda'\otimes_{R'}(\Lambda')\op}
\Lambda'\cong R'$ and $R'$ is local, there exists some $i,j$ such that
$u_iv_j\psi_i\varphi_j$ is an automorphism of $\Lambda'$.
Then $\psi_i:M'\rightarrow \Lambda'$ is also an isomorphism.
By faithful flatness, $\psi_i: M\rightarrow \Lambda$ is an isomorphism.

{\bf 2}.
It is easy to see that $M\in\md\Lambda$ is projective.
So replacing $\Lambda$ by $\Lambda/J$, where $J$ is the radical of $J$,
and changing $R$ and $R'$ as above, we may assume that $R$ is a field and
$\Lambda$ is central simple.
Then there is only one simple right $\Lambda$-module,
and $M$ and $\Lambda$ are direct sums of copies of it.
As $M'\cong\Lambda'$,
by dimension counting, the number of copies are equal,
and hence $M$ and $\Lambda$ are isomorphic.
\end{proof}

\begin{lemma}\label{2-canonical-descent.lem}
  Let $(R,\fm,k)\rightarrow (R',\fm',k')$ be a flat local homomorphism
  between Noetherian local rings.
  \begin{enumerate}
  \item[\bf 1] Let $C$ be a $2$-canonical bimodule of $\Lambda$ over $R$.
    Let $M$ be a $\Lambda$-bimodule such that $M':=R'\otimes_R M$
    is isomorphic to $C':=R'\otimes_R C$ as a $\Lambda'$-bimodule.
    Then $M\cong C$ as a $\Lambda$-bimodule.
  \item[\bf 2] Let $C$ be a $2$-canonical right
    $\Lambda$-module over $R$.
    Let $M$ be a right $\Lambda$-module such that $M':=R'\otimes_R M$
    is isomorphic to $C':=R'\otimes_R C$ as a right $\Lambda'$-module.
    Then $M\cong C$ as a right $\Lambda$-module.
  \end{enumerate}
\end{lemma}

\begin{proof}
  {\bf 1}.
  As $M'\cong C'$ and $C\in (S_2')_C$, it is easy to see that $M\in (S_2')_C$.
  Hence $M$ is a $\Gamma$-bimodule, where $\Gamma=\End\Lop C=\End_{\Lambda}C$, see
  (\ref{Gamma-Delta.par}) and Corollary~\ref{2-canonical-equiv.cor}.
  Note that $(M\da)' \cong (C\da)'\cong \Gamma'$ as $\Gamma'$-bimodules.
  By Lemma~\ref{flat-local-free.lem}, {\bf 1}, we have that $M\da\cong\Gamma$
  as a $\Gamma$-bimodule.
  Hence $M\cong M\ddd \cong \Gamma\dda\cong C$.

  {\bf 2}.
  As $(M\da)'\cong (C\da)'\cong \Gamma'$ as $\Gamma'$-modules,
  $M\da\cong \Gamma$ as $\Gamma$-modules by
  Lemma~\ref{flat-local-free.lem}, {\bf 2}.
  Hence $M\cong M\ddd\cong \Gamma\dda\cong C$.
\end{proof}

\begin{proposition}\label{canonical-descent.prop}
  Let $(R,\fm,k)\rightarrow (R',\fm',k')$
  be a flat local homomorphism between Noetherian local rings.
  Assume that $R'/\fm R'$ is zero-dimensional, and
  $M':=R'\otimes_R M$ is the right canonical module
  of $\Lambda':=R'\otimes_R \Lambda$.
  Then $R'/\fm R'$ is Gorenstein.
\end{proposition}

\begin{proof}
  We may assume that both $R$ and $R'$ are complete.
  Replacing $R$ by $R/\ann_R \Lambda$ and $R'$ by $R'\otimes_R R/\ann_R
  \Lambda$, we may assume that $\Lambda$
  is a faithful $R$-module.
  Let $d=\dim R=\dim R'$.

  Then
  \[
  R'\otimes_R H^d_\fm(M)\cong H^d_{\fm'}(R'\otimes_R M)\cong H^d_{\fm'}(K_{\Lambda'})
  \cong \Hom_{R'}(\Gamma',E'),
  \]
  where $\Lambda'=R'\otimes_R \Lambda$,
  $E'=E_{R'}(R'/\fm')$ is the injective hull of the residue field,
  $\Gamma=\End\Lop M$, $\Gamma'=R'\otimes_R \Gamma\cong \End_{\Lambda'}K_{\Lambda'}$,
  and the isomorphisms are those of $\Gamma'$-modules.
  The last isomorphism is by (\ref{canonical-top-local-cohomology.par}).
  So $R'\otimes_R H^d_\fm(M)\in \Mod \Gamma'$ is injective.
  Considering the spectral sequence
  \begin{multline*}
    E_2^{p,q}=\Ext_{R'\otimes_R(\Gamma\otimes_R k)}^p(W,
                 \Ext^q_{\Gamma'}(R'\otimes_R (\Gamma\otimes_R k),
                                         R'\otimes_R H^d_\fm(M)))\\
  \Rightarrow
  \Ext^{p+q}_{\Gamma'}(W,R'\otimes_R H^d_\fm(M))
  \end{multline*}
  for $W\in \Mod (R'\otimes_R (\Gamma\otimes_R k))$, 
  $E_2^{1,0}=E_\infty^{1,0}\subset \Ext^1_{\Gamma'}(W,R'\otimes_R H^d_\fm(M))=0$
  by the injectivity of $R'\otimes_R H^d_\fm(M)$.
  It follows that
  \[
  \Hom_{\Gamma'}(R'\otimes_R (\Gamma\otimes_R k),R'\otimes_R H^d_\fm(M))
  \cong
  (R'/\fm R')\otimes_k \Hom_R(k,H^d_\fm(M))
  \]
  is an injective $(R'/\fm R') \otimes_k (\Gamma\otimes_R k)$-module.
  However, as an $R'/\fm R'$-module, this is a free module.
  Also, this module must be an injective $R'/\fm R'$-module, and hence
  $R'/\fm R'$ must be Gorenstein.
\end{proof}

\begin{lemma}\label{canonical-ascent.lem}
  Let $(R,\fm,k)\rightarrow (R',\fm',k')$ be a flat local homomorphism
  between Noetherian local rings such that $R'/\fm R'$ is Gorenstein.
  Assume that the canonical module $K_\Lambda$ of $\Lambda$ exists.
  Then $R'\otimes_R K_\Lambda$ is the canonical module of $R'\otimes_R \Lambda$.
\end{lemma}

\begin{proof}
  We may assume that both $R$ and $R'$ are complete.
  Let $\Bbb I$ be the normalized dualizing complex of $R$.
  Then $R'\otimes_R \Bbb I[d'-d]$ is a normalized
  dualizing complex of $R'$,
  where $d'=\dim R'$ and $d=\dim R$,
  since $R\rightarrow R'$ is a flat local homomorphism with the $d'-d$-dimensional
  Gorenstein closed fiber, see \cite[(5.1)]{AF} (the definition of a normalized
  dualizing complex in \cite{AF} is different from ours. We follow
  the one in \cite[Chapter~V]{Hartshorne}).
  So
  \[
  R'\otimes_R K_\Lambda \cong R'\otimes_R \Ext^{-d}_R(\Lambda,\Bbb I)
  \cong \Ext^{-d'}_R(R'\otimes_R \Lambda,R'\otimes_R \Bbb I[d'-d])
  \cong K_{\Lambda'}.
  \]
\end{proof}

\begin{theorem}[(Non-commutative Aoyama's theorem) cf.~{\cite[Theorem~4.2]{Aoyama}}]
  \label{Aoyama.thm}
  Let $(R,\fm)\rightarrow (R',\fm')$ be a flat local homomorphism between
  Noetherian local rings.
  \begin{enumerate}
  \item[\bf 1] If $M$ is a $\Lambda$-bimodule and $M'=R'\otimes_R M$ is
    the canonical module of $\Lambda'=R'\otimes_R \Lambda$, then
    $M$ is the canonical module of $\Lambda$.
  \item[\bf 2] If $M$ is a right $\Lambda$-module such that $M'$ is
    the right canonical module of $\Lambda'$, then $M$ is the right
    canonical module of $\Lambda$.
  \end{enumerate}
\end{theorem}

\begin{proof}
  We may assume that both $R$ and $R'$ are complete.
  Then the canonical module exists, and the
  localization of a canonical module is a canonical module, and hence
  we may localize $R'$ by a minimal element of
  $\{P\in \Spec R'\mid P\cap R=\fm\}$, and take the completion again, we
  may further assume that the fiber ring $R'/\fm R'$ is zero-dimensional.
  Then $R'/\fm R'$ is Gorenstein by Proposition~\ref{canonical-descent.prop}.
  Then by Lemma~\ref{canonical-ascent.lem},
  $M'\cong K_{\Lambda'}\cong R'\otimes_R K_\Lambda$.
  By Lemma~\ref{2-canonical-descent.lem}, $M\cong K_\Lambda$.
  In {\bf 1}, the isomorphisms are those of bimodules, while in {\bf 2},
  they are of right modules.
  The proofs of {\bf 1} and {\bf 2} are complete.
\end{proof}

\begin{corollary}\label{canonical-localization.cor}
  Let $(R,\fm)$ be a Noetherian local ring, and assume that $K$ is the
  canonical \(resp.\ right canonical\) module of $\Lambda$.
  If $P\in \Supp_R K$, then the localization $K_P$ is the
  canonical \(resp.\ right canonical\) module of $\Lambda_P$.
  In particular, $K$ is a semicanonical bimodule \(resp.\ right module\),
  and hence is $2$-canonical over $R/\ann_R\Lambda$.
\end{corollary}

\begin{proof}
  Let $Q$ be a prime ideal of $\hat R$ lying over $P$.
  Then $(\hat K)_Q\cong \hat R_Q  \otimes_{R_P} K_P$ is nonzero by
  assumption, and hence is the
  canonical (resp.\ right canonical)
  module of $\hat R_Q \otimes_R \Lambda$.
  Using Theorem~\ref{Aoyama.thm}, $K_P$ is
  the canonical (resp.\ right canonical) module of $\Lambda_P$.
  The last assertion follows.
\end{proof}

\paragraph
Let $(R,\fm)$ be local, and assume that $K_\Lambda$ exists.
Assume that $\Lambda$ is a faithful $R$-module.
Then it is a $2$-canonical $\Lambda$-bimodule over $R$ by
Corollary~\ref{canonical-localization.cor}.
Letting $\Gamma=\End\Lop K_\Lambda$, 
$K_\Gamma \cong K_\Lambda$ as $\Lambda$-bimodules by
Corollary\ref{lambda-gamma2.cor}.
So by Corollary~\ref{2-canonical-equiv.cor},
there exists some $\Gamma$-bimodule structure of $K_\Lambda$ such
that $K_\Gamma \cong K_\Lambda$ as $\Gamma$-bimodules.
As the left $\Gamma$-module structure of $K_\Lambda$ which
extends the original left $\Lambda$-module structure is unique,
and it is the obvious action of $\Gamma=\End\Lop K_\Lambda$.
Similarly the right action of $\Gamma$ is the obvious
action of $\Gamma=\Delta=(\End_{\Lambda}K_\Lambda)\op$,
see (\ref{Gamma-Delta.par}).

\section{Evans--Griffith's theorem for $n$-canonical modules}\label{Evans-Griffith.sec}

\begin{lemma}[cf.~{\cite[Proposition~2]{Aoyama}, \cite[Proposition~4.2]{Ogoma},
  \cite[Proposition~1.2]{AG}}]\label{AG.lem}
  Let $(R,\fm)$ be local and assume that $\Lambda$ has a canonical module
  $C=K_\Lambda$.
  Then we have
  \begin{enumerate}
  \item[\bf 1] $\lambda_R:\Lambda\rightarrow \End\Lop K_\Lambda$ is injective
    if and only if $\Lambda$ satisfies the $(S_1)^R$ condition and 
    $\Supp_R\Lambda$ is equidimensional.
  \item[\bf 2] $\lambda_R:\Lambda\rightarrow \End\Lop K_\Lambda$ is bijective
    if and ony if $\Lambda$ satisfies the $(S_2)^R$ condition.
  \end{enumerate}
\end{lemma}

\begin{proof}
  Replacing $R$ by $R/\ann_R \Lambda$, we may assume that $\Lambda$ is a
  faithful $R$-module.
  Then $K_\Lambda$ is a $2$-canonical $\Lambda$-bimodule over $R$ by
  Corollary~\ref{canonical-localization.cor}.
  $K_\Lambda$ is full if and only if $\Supp_R\Lambda$ is equidimensional
  by Lemma~\ref{non-complete.lem}, {\bf 1}.

  Now {\bf 1} is a consequence of Lemma~\ref{1-canonical-injective.lem}.
  {\bf 2} follows from Corollary~\ref{2-canonical-isom.cor} and
  Lemma~\ref{ogoma.lem}.
\end{proof}

\begin{proposition}[cf.~{\cite[(2.3)]{AG}}]\label{AG.prop}
  Let $(R,\fm)$ be a local ring, and assume that there is an $R$-canonical module
  $K_\Lambda$ of $\Lambda$.
  Assume that $\Lambda\in (S_2)^R$, and $K_\Lambda$ is a Cohen--Macaulay $R$-module.
  Then $\Lambda$ is Cohen--Macaulay.
  If, moreover, $K_\Lambda$ is maximal Cohen--Macaulay, then so is $\Lambda$.
\end{proposition}
  
\begin{proof}
  The second assertion follows from the first.
  We prove the first assertion.
  Replacing $R$ by $R/\ann_R\Lambda$, we may assume that $\Lambda$ is faithful.
  Let $d=\dim R$.
  So $\Lambda$ satisfies $(S_2')$, and $K_\Lambda$ is maximal Cohen--Macaulay.
  As $K_\Lambda$ is the lowest non-vanishing cohomology of
  $\Bbb J:=\RHom_R(\Lambda,\Bbb I)$, there is a natural map
  $\sigma:K_\Lambda[d]\rightarrow \Bbb J$ which induces an isomorphism on the
  $-d$th cohomology groups.
  Then the diagram
  \[
  \xymatrix{
    \Lambda \ar[r]^-{\lambda} \ar[d]^{\lambda} &
    \Hom_{\Lambda\op}(K_\Lambda[d],K_\Lambda[d]) \ar[d]^{\sigma_*} \\
    \RHom\Lop(\Bbb J,\Bbb J) \ar[r]^-{\sigma^*} &
    \RHom\Lop(K_\Lambda[d],\Bbb J)
  }
  \]
  is commutative.
  The top horizontal arrow $\lambda$ is an isomorphism by Lemma~\ref{AG.lem}.
  Note that
  \[
  \RHom\Lop(\Bbb J,\Bbb J)\cong \RHom_R(\Bbb J,\Bbb I)
  =\RHom_R(\RHom_R(\Lambda,\Bbb I),\Bbb I)=\Lambda,
  \]
  and the left vertical arrow is an isomorphism.
  As $K_\Lambda$ is maximal Cohen--Macaulay, $\RHom\Lop(K_\Lambda[d],\Bbb J)$ is
  concentrated in degree zero.
  As $H^{i}(\Bbb J)=0$ for $i<-d$, we have that the right vertical arrow
  $\sigma_*$ is an isomorphism.
  Thus the bottom horizontal arrow $\sigma^*$ is an isomorphism.
  Applying $\RHom_\Lambda(?,\Bbb J)$ to this map, we have that
  $K_\Lambda[d]\rightarrow \Bbb J$ is an isomorphism.
  So $\Lambda$ is Cohen--Macaulay, as desired.
  \end{proof}

\begin{corollary}[cf.~{\cite[(2.2)]{AG}}]\label{AG.cor}
  Let $(R,\fm)$ be a local ring, and assume that there is an $R$-canonical module
  $K_\Lambda$ of $\Lambda$.
  Then $K_\Lambda$ is a Cohen--Macaulay \(resp.\ maximal Cohen--Macaulay\)
  $R$-module if and only if $\Gamma=\End_{\Lambda\op}K_\Lambda$ is so.
\end{corollary}

\begin{proof}
  As $K_\Lambda$ and $\Gamma$ has the same support, if both of them are Cohen--Macaulay
  and one of them are maximal Cohen--Macaulay, then the other is also.
  So it suffices to prove the assertion on the Cohen--Macaulay property.
  To verify this, we may assume that $\Lambda$ is a faithful $R$-module.
  Note that $\Gamma$ satisfies $(S_2')$.
  By Corollary~\ref{lambda-gamma2.cor}, $K_\Lambda$ is Cohen--Macaulay if and
  only if $K_\Gamma$ is.
  If $\Gamma$ is Cohen--Macaulay, then $K_\Gamma$ is Cohen--Macaulay by
  (\ref{algebra-isom.par}).
  Conversely, if $K_\Gamma$ is Cohen--Macaulay, then $\Gamma$ is Cohen--Macaulay
  by Proposition~\ref{AG.prop}.
\end{proof}

\begin{theorem}[cf.~{\cite[(3.8)]{EG}, \cite[(3.1)]{AI}}]\label{main.thm}
  Let $R$ be a Noetherian commutative ring, and $\Lambda$ a module-finite
  $R$-algebra, which may not be commutative.
  Let $n\geq 1$, and $C$ be a right $n$-canonical $\Lambda$-module.
  Set $\Gamma=\End\Lop C$.
  Let $M\in\md C$.
  Then the following are equivalent.
  \begin{enumerate}
  \item[\bf 1] $M\in \TF(n,C)$.
  \item[\bf 2] $M\in \UP(n,C)$.
  \item[\bf 3] $M\in\Syz(n,C)$.
  \item[\bf 4] $M\in (S_n')_C$.
  \end{enumerate}
\end{theorem}

\begin{proof}
  {\bf 1$\Rightarrow$2$\Rightarrow$3$\Rightarrow$4} is easy.
  We prove {\bf 4$\Rightarrow$1}.
  By Lemma~\ref{1-canonical.lem}, we may assume that $n\geq 2$.
  By Lemma~\ref{2-canonical.lem}, $M\in\TF(2,C)$.
  Let
  \[
  \Bbb F:
  0\leftarrow M\da \leftarrow F_0 \leftarrow F_1 \leftarrow\cdots\leftarrow F_{n-1}
  \]
  be a resolution of $M\da$ in $\Gamma\md$ with each $F_i\in\add\Gamma$.
  It suffices to prove its dual
  \[
  \Bbb F\dda:
  0\rightarrow M \rightarrow F_0\dda \rightarrow F_1\dda \rightarrow\cdots\rightarrow
  F_{n-1}\dda
  \]
  is acyclic.
  By Lemma~\ref{s_n-acyclicity.lem}, we may localize at $P\in R\an{<n}$, and
  may assume that $\dim R<n$.
  If $M=0$, then $\Bbb F$ is split exact, and so $\Bbb F\dda$ is also exact.
  So we may assume that $M\neq 0$.
  Then by assumption, $C\cong K_\Lambda$ in $\md \Lambda$, and $C$ is
  a maximal Cohen--Macaulay $R$-module.
  Hence $\Gamma$ is Cohen--Macaulay by Corollary~\ref{AG.cor}.
  So by (\ref{ext-vanishing.par}) and Lemma~\ref{lambda-gamma2.cor}, 
  $\RHom_\Gamma(M\da,C)=\RHom_\Gamma(M\da,K_\Gamma)=M$, and we are done.
\end{proof}
  
\begin{corollary}\label{main.cor}
  Let the assumptions and notation be as in {\rm Theorem~\ref{main.thm}}.
  Let $n\geq 0$.
  Assume further that
  \begin{enumerate}
  \item[\bf 1] $\Ext^i\Lop(C,C)=0$ for $1\leq i\leq n$;
  \item[\bf 2] $C$ is $\Lambda$-full.
  \item[\bf 3] $\Lambda$ satisfies the $(S_n')^R$ condition.
  \end{enumerate}
  Then for $0\leq r\leq n$, $\sph r C$ is contravariantly finite
  in $\md\Lambda$.
\end{corollary}

\begin{proof}
  For any $M\in\md \Lambda$, the $n$th syzygy module
  $\Omega^n M$ satisfies the $(S_n')^R_C$-condition by {\bf 2} and {\bf 3}.
  By Theorem~\ref{main.thm}, $\Omega^n M\in \TF\Lop(n,C)$.
  By Theorem~\ref{X_{m,n}-approximation.thm}, $M\in\Z_{r,0}$, and there
  is a short exact sequence
  \[
  0\rightarrow Y\rightarrow X\xrightarrow g M\rightarrow 0
  \]
  with $X\in \X_{r,0}=\sph r C$ and $Y\in \Y_r$.
  As $\Ext^1\Lop(X,Y)=0$, we have that $g$ is a right $\sph r C$-approximation,
  and hence $\sph r C$ is contravariantly finite.
\end{proof}

\begin{corollary}\label{main2.cor}
  Let the assumptions and notation be as in {\rm Theorem~\ref{main.thm}}.
  Let $n\geq 0$, and $C$ a $\Lambda$-full $(n+2)$-canonical $\Lambda$-bimodule
  over $R$.
  Assume that $\Lambda$ satisfies the $(S_{n+2}')^R$ condition.
  Then $\sph n C$ is contravariantly finite in $\md \Lambda$.
\end{corollary}

\begin{proof}
  By Corollary~\ref{main.cor}, it suffices to show that $\Ext^i\Lop(C,C)=0$
  for $1\leq i\leq n$.
  Let $\Delta=(\End_\Lambda C)\op$.
  Then the canonical map $\Lambda\rightarrow \Delta$ is an isomorphism
  by Lemma~\ref{2-canonical-isom.cor}, 
  since $C$ is a $\Lambda$-full $2$-canonical $\Lambda$-bimodule over $R$.
  As $\Lambda\in (S_{n+2}')^R$ and $C$ is a $\Lambda$-full $(n+2)$-canonical
  left $\Lambda$-module over $R$, applying Theorem~\ref{main.thm} to $\Lambda\op$,
  we have that $\Ext^i_{\Delta\op}(C,C)=0$ for $1\leq i\leq n$.
  As we have $\Lambda\op\rightarrow \Delta\op$ is an isomorphism,
  we have that $\Ext^i\Lop(C,C)=0$, as desired.
\end{proof}

\section{Symmetric and Frobenius algebras}\label{Frobenius.sec}

\paragraph Let $(R,\fm)$ be a Noetherian semilocal ring, and $\Lambda$ a
module-finite $R$-algebra.
We say that $\Lambda$ is {\em quasi-symmetric}
if $\Lambda$ is the canonical module of $\Lambda$.
That is, $\Lambda\cong K_\Lambda$ as $\Lambda$-bimodules.
It is called {\em symmetric} if it is quasi-symmetric and GCM.
Note that $\Lambda$ is quasi-symmetric (resp.\ symmetric) if and only
if $\hat \Lambda$ is so, where $\hat ?$ denotes the $\fm$-adic completion.
Note also that quasi-symmetric and symmetric are absolute notion, and
is independent of the choice of $R$ in the sense that the definition
does not change when we replace $R$ by the center of $\Lambda$.

\paragraph
For (non-semilocal) Noetherian ring $R$,
we say that $\Lambda$ is locally quasi-symmetric
(resp.\ locally symmetric) over $R$ if for any $P\in \Spec R$,
$\Lambda_P$ is a quasi-symmetric (resp.\ symmetric) $R_P$-algebra.
This is equivalent to say that for any maximal ideal $\fm $ of $R$,
$\Lambda_\fm$ is quasi-symmetric (resp.\ symmetric).
In the case that $(R,\fm)$ is semilocal,
$\Lambda$ is 
locally quasi-symmetric (resp.\ locally symmetric) over $R$ if
it is quasi-symmetric (resp.\ symmetric), but the converse is
not true in general.

\begin{lemma}\label{pseudo-Frobenius.lem}
  Let $(R,\fm)$ be a Noetherian
  semilocal ring, and $\Lambda$ a module-finite
  $R$-algebra.
  Then the following are equivalent.
  \begin{enumerate}
  \item[\bf 1] $\Lambda_\Lambda$ is the right canonical module of $\Lambda$.
  \item[\bf 2] ${}_\Lambda\Lambda$ is the left canonical module of $\Lambda$.
  \end{enumerate}
\end{lemma}

\begin{proof}
  We may assume that $R$ is complete.
  Then replacing $R$ by a Noether normalization of
  $R/\ann_R \Lambda$, we may assume
  that $R$ is regular and $\Lambda$ is a faithful $R$-module.

  We prove {\bf 1$\Rightarrow$\bf 2}.
  By Lemma~\ref{non-complete.lem}, $\Lambda$ satisfies $(S_2')^R$.
  As $R$ is regular and $\dim R=\dim \Lambda$, $K_\Lambda=\Lambda^*
  =\Hom_R(\Lambda,R)$.
  So we get an $R$-linear map
  \[
  \varphi:\Lambda\otimes_R\Lambda \rightarrow R
  \]
  such that $\varphi(ab\otimes c)=\varphi(a\otimes bc)$ and that
  the induced map $h:\Lambda\rightarrow \Lambda^*$ given by $h(a)(c)=\varphi(a\otimes
  c)$ is an isomorphism (in $\md \Lambda$).
  Now $\varphi$ induces a homomorphism $h':\Lambda\rightarrow \Lambda^*$ in
  $\Lambda\md$ given by $h'(c)(a)=\varphi(a\otimes c)$.
  To verify that this is an isomorphism, as $\Lambda$ and $\Lambda^*$ are
  reflexive $R$-modules, we may localize at $P\in R\an{<2}$, and then take a
  completion, and hence we may further assume that $\dim R\leq 1$.
  Then $\Lambda$ is a finite free $R$-module, and the matrices of $h$ and $h'$ are
  transpose each other.
  As the matrix of $h$ is invertible, so is that of $h'$, and $h'$ is an isomorphism.

  {\bf 2$\Rightarrow$1} follows from {\bf 1$\Rightarrow$2},
  considering the opposite ring.
\end{proof}

\begin{definition}
  Let $(R,\fm)$ be semilocal.
  We say that $\Lambda$ is a {\em pseudo-Frobenius $R$-algebra}
  if the equivalent conditions of Lemma~\ref{pseudo-Frobenius.lem} are
  satisfied.
  If $\Lambda$ is GCM in addition, then it is called a
  {\em Frobenius $R$-algebra}.
  Note that these definitions
  are independent of the choice of $R$.
  Moreover, $\Lambda$ is pseudo-Frobenius (resp.\ Frobenius) if and only
  if $\hat\Lambda$ is so, where $\hat ?$ is the $\fm$-adic completion.
  For a general $R$, we say that $\Lambda$ is locally pseudo-Frobenius
  (resp.\ locally Frobenius) over $R$
  if $\Lambda_P$ is pseudo-Frobenius (resp.\ Frobenius) for $P\in\Spec R$.
\end{definition}

\begin{lemma}\label{quasi-Frobenius.lem}
  Let $(R,\fm)$ be semilocal.
  Then the following are equivalent.
  \begin{enumerate}
  \item[\bf 1] 
    $(K_{\hat\Lambda})_{\hat\Lambda}$ is projective in $\md\hat\Lambda$.
  \item[\bf 2] 
    ${}_{\hat\Lambda} (K_{\hat \Lambda})$ 
    is projective in $\hat\Lambda\md$,
  \end{enumerate}
  where $\hat ?$ denotes the $\fm$-adic completion.
\end{lemma}

\begin{proof}
  We may assume that $(R,\fm,k)$ is complete regular local and
  $\Lambda$ is a faithful $R$-module.
  Let $\bar ?$ denote the functor $k\otimes_R ?$.
  Then $\bar \Lambda$ is a finite dimensional $k$-algebra.
  So $\md\bar\Lambda$ and $\bar\Lambda\md$ have the same number of
  simple modules, say $n$.
  An indecomposable projective module in $\md\Lambda$ is nothing but
  the projective cover of a simple module in $\md\bar\Lambda$.
  So $\md\Lambda$ and $\Lambda\md$ have $n$ indecomposable projectives.
  Now $\Hom_R(?,R)$ is an equivalence between $\add (K_\Lambda)_\Lambda$ and
  $\add {}_\Lambda \Lambda$.
  It is also an equivalence between $\add {}_\Lambda (K_\Lambda)$ and
  $\add\Lambda_\Lambda$.
  So both $\add(K_\Lambda)_\Lambda$ and $\add{}_\Lambda (K_\Lambda)$ also have
  $n$ indecomposables.
  So {\bf 1} is equivalent to $\add (K_\Lambda)_\Lambda = \add \Lambda_\Lambda$.
  {\bf 2} is equivalent to $\add {}_\Lambda (K_\Lambda) = \add {}_\Lambda\Lambda$.
  So {\bf 1$\Leftrightarrow $2} is proved simply applying
  the duality $\Hom_R(?,R)$.
\end{proof}

\paragraph
Let $(R,\fm)$ be semilocal.
If the equivalent conditions in Lemma~\ref{quasi-Frobenius.lem} are
satisfied, then we say that $\Lambda$ is {\em pseudo-quasi-Frobenius}.
If it is GCM in addition, then we say that it is
{\em quasi-Frobenius}.
These definitions are independent of the choice of $R$.
Note that $\Lambda$ is pseudo-quasi-Frobenius (resp.\ quasi-Frobenius)
if and only if $\hat\Lambda$ is so.

\begin{proposition}\label{GN.prop}
  Let $(R,\fm)$ be semilocal.
  Then the following are equivalent.
  \begin{enumerate}
  \item[\bf 1] $\Lambda$ is quasi-Frobenius.
  \item[\bf 2] $\Lambda$ is GCM, and $\dim\Lambda=\idim {}_\Lambda\Lambda$,
    where $\idim$ denotes the injective dimension.
  \item[\bf 3] $\Lambda$ is GCM, and $\dim\Lambda=\idim \Lambda_\Lambda$.
  \end{enumerate}
\end{proposition}

\begin{proof}
  {\bf 1$\Rightarrow $2}.
  By definition, $\Lambda$ is GCM.
  To prove that $\dim\Lambda=\idim {}_\Lambda\Lambda$, we may assume that
  $R$ is local.
  Then by \cite[(3.5)]{GN}, we may assume that $R$ is complete.
  Replacing $R$ by the Noetherian normalization of
  $R/\ann_R \Lambda$, we may
  assume that $R$ is a complete regular local ring of dimension $d$, and
  $\Lambda$ its maximal Cohen--Macaulay module.
  As $\add {}_\Lambda \Lambda = \add {}_\Lambda (K_\Lambda)$ by
  the proof of Lemma~\ref{quasi-Frobenius.lem}, it suffices to
  prove $\idim {}_\Lambda (K_\Lambda)=d$.
  Let $\Bbb I_R$ be the minimal injective resolution of the
  $R$-module $R$.
  Then $\Bbb J=\Hom_R(\Lambda,\Bbb I_R)$
  is an injective resolution
  of $K_\Lambda=\Hom_R(\Lambda,R)$.
  As the length of $\Bbb J$ is $d$ and
  \[
  \Ext^d_\Lambda(\Lambda/\fm\Lambda,K_\Lambda)\cong
  \Ext^d_R(\Lambda/\fm\Lambda,R)\neq 0,
  \]
  we have that $\idim{}_\Lambda (K_\Lambda)=d$.

  {\bf 2$\Rightarrow$1}.
  We may assume that $R$ is complete regular local and
  $\Lambda$ is maximal Cohen--Macaulay.
  By \cite[(3.6)]{GN}, we may further assume that $R$ is a field.
  Then ${}_\Lambda \Lambda$ is injective.
  So $(K_\Lambda)_\Lambda=\Hom_R(\Lambda,R)$ is projective,
  and $\Lambda$ is quasi-Frobenius, see \cite[(IV.3.7)]{SY}.

  {\bf 1$\Leftrightarrow$3} is proved similarly.
\end{proof}

\begin{corollary}\label{GN.cor}
  Let $R$ be arbitrary.
  Then the following are equivalent.
  \begin{enumerate}
  \item[\bf 1] For any $P\in\Spec R$, $\Lambda_P$ is quasi-Frobenius.
  \item[\bf 2] For any maximal ideal $\fm$ of $R$, $\Lambda_\fm$ is quasi-Frobenius.
  \item[\bf 3] $\Lambda$ is a Gorenstein $R$-algebra in the sense that
    $\Lambda$ is a Cohen--Macaulay $R$-module, and
    $\idim_{\Lambda_P} {}_{\Lambda_P}\Lambda_P=\dim \Lambda_P$ for
    any $P\in\Spec R$.
  \end{enumerate}
\end{corollary}

\begin{proof}
  {\bf 1$\Rightarrow$2} is trivial.

  {\bf 2$\Rightarrow$3}.
  By Proposition~\ref{GN.prop}, we have
  $\idim {}_{\Lambda_\fm}\Lambda_\fm=\dim \Lambda_\fm$
  for each $\fm$.
  Then by \cite[(4.7)]{GN}, $\Lambda$ is a Gorenstein $R$-algebra.

  {\bf 3$\Rightarrow$1} follows from Proposition~\ref{GN.prop}.
\end{proof}
  
\paragraph
Let $R$ be arbitrary.
We say that $\Lambda$ is a {\em quasi-Gorenstein} $R$-algebra if
$\Lambda_P$ is pseudo-quasi-Frobenius for each $P\in\Spec R$.

\begin{definition}[Scheja--Storch \cite{SS}]
Let $R$ be general.
We say that $\Lambda$ is symmetric (resp.\ Frobenius)
relative to
$R$ if $\Lambda$ is $R$-projective, and $\Lambda^*:=\Hom_R(\Lambda,R)$ is
isomorphic to $\Lambda$ as a $\Lambda$-bimodule (resp.\ as a right $\Lambda$-module).
It is called quasi-Frobenius relative to $R$
if the right $\Lambda$-module $\Lambda^*$ is projective.
\end{definition}

\begin{lemma}\label{SS.lem}
  Let $(R,\fm)$ be local.
  \begin{enumerate}
  \item[\bf 1] If $\dim \Lambda=\dim R$, $R$ is quasi-Gorenstein, and
    $\Lambda^* \cong \Lambda$ as $\Lambda$-bimodules
    \(resp.\ $\Lambda^*\cong\Lambda$ as right $\Lambda$-modules,
    $\Lambda^*$ is projective as a right $\Lambda$-module\), then
    $\Lambda$ is quasi-symmetric \(resp.\ pseudo-Frobenius, pseudo-quasi-Frobenius\).
  \item[\bf 2] If $R$ is Gorenstein and $\Lambda$ is symmetric \(resp.\
    Frobenius, quasi-Frobenius\) relative to $R$, then $\Lambda$ is symmetric
    \(resp.\ Frobenius, quasi-Frobenius\).
  \item[\bf 3] If $\Lambda$ is nonzero and $R$-projective, then
    $\Lambda$ is quasi-symmetric \(resp.\ pseudo-Frobenius,
    pseudo-quasi-Frobenius\) if and only if $R$ is quasi-Gorenstein and
    $\Lambda$ is symmetric \(resp.\ Frobenius, quasi-Frobenius\)
    relative to $R$.
  \item[\bf 4] If $\Lambda$ is nonzero and $R$-projective, then
    $\Lambda$ is symmetric \(resp.\ Frobenius, quasi-Frobenius\) if and only if
    $R$ is Gorenstein and $\Lambda$ is symmetric \(resp.\ Frobenius, quasi-Frobenius\)
    relative to $R$.
  \end{enumerate}
\end{lemma}

\begin{proof}
  We can take the completion, and we may assume that $R$ is complete local.
  
  {\bf 1}.
  Let $d=\dim \Lambda=\dim R$, and let $\Bbb I$ be the normalized dualizing complex of $R$.
  Then
  \[
  K_\Lambda= \Ext^{-d}_R(\Lambda,\Bbb I)\cong \Hom_R(\Lambda,H^{-d}(\Bbb I))
  \cong \Hom(\Lambda,K_R)\cong \Hom(\Lambda,R)=\Lambda^*
  \]
  as $\Lambda$-bimodules, and the result follows.

  {\bf 2}.
  We may assume that $\Lambda$ is nonzero.
  As $R$ is Cohen--Macaulay and $\Lambda$ is a finite projective $R$-module,
  $\Lambda$ is a maximal Cohen--Macaulay $R$-module.
  By {\bf 1}, the result follows.

  {\bf 3}.
  The \lq if' part follows from {\bf 1}.
  We prove the \lq only if' part.
  As $\Lambda$ is $R$-projective and nonzero, $\dim\Lambda=\dim R$.
  As $\Lambda$ is $R$-finite free, $K_\Lambda\cong \Hom_R(\Lambda,K_R)\cong
  \Lambda^*\otimes_R K_R$.
  As $K_\Lambda$ is $R$-free and $\Lambda^*\otimes_R K_R$ is nonzero and is isomorphic to
  a direct sum of copies of $K_R$, we have that $K_R$ is $R$-projective,
  and hence $R$ is quasi-Gorenstein, and $K_R\cong R$.
  Hence $K_\Lambda\cong \Lambda^*$, and the result follows.

  {\bf 4} follows from {\bf 3} easily.
\end{proof}

\paragraph
Let $(R,\fm)$ be semilocal.
Let a finite group $G$ act on $\Lambda$ by $R$-algebra automorphisms.
Let $\Omega=\Lambda*G$, the twisted group algebra.
That is, $\Omega=\Lambda\otimes_R RG=\bigoplus_{g\in G}\Lambda g$
as an $R$-module, and
the product of $\Omega$ is given by $(ag)(a'g')=(a(ga'))(gg')$
for $a,a'\in\Lambda$ and $g,g'\in G$.
This makes $\Omega$ a module-finite $R$-algebra.

\paragraph
We simply call an $RG$-module a $G$-module.
We say that $M$ is a $(G,\Lambda)$-module if $M$ is a $G$-module,
$\Lambda$-module, the $R$-module structures coming from that of the
$G$-module structure and the $\Lambda$-module structure agree, and
$g(am)=(ga)(gm)$ for $g\in G$, $a\in \Lambda$, and $m\in M$.
A $(G,\Lambda)$-module and an $\Omega$-module are one and the same thing.

\paragraph
By the action $(a\otimes a')g)a_1=a(ga_1)a'$,
we have that $\Lambda$ is a $(\Lambda\otimes\Lambda\op)*G$-module
in a natural way.
So it is an $\Omega$-module by the action $(ag)a_1=a(ga_1)$.
It is also a right $\Omega$-module by the action $a_1(ag)=g^{-1}(a_1a)$.
If the action of $G$ on $\Lambda$ is trivial, then these actions
make an $\Omega$-bimodule.

\paragraph
Given an $\Omega$-module $M$ and an $RG$-module $V$,
$M\otimes_R V$ is an $\Omega$-module by
$(ag)(m\otimes v)=(ag)m\otimes gv$.
$\Hom_R(M,V)$ is a right $\Omega$-module by
$(\varphi(ag))(m)= g^{-1}\varphi(a(gm))$.
It is easy to see that the standard isomorphism
\[
\Hom_R(M\otimes_R V,W)\rightarrow \Hom_R(M,\Hom_R(V,W))
\]
is an isomorphism of right $\Omega$-modules for a left $\Omega$-module
$M$ and $G$-modules $V$ and $W$.

\paragraph
Now consider the case $\Lambda=R$.
Then the pairing $\phi: RG\otimes_R RG \rightarrow R$ given by
$\phi(g\otimes g')=\delta_{gg',e}$ (Kronecker's delta)
is non-degenerate, and induces
an $RG$-bimodule isomorphism $\Omega=RG\rightarrow (RG)^*=\Omega^*$.
As $\Omega=RG$ is a finite free $R$-module, we have that
$\Omega=RG$ is symmetric relative to $R$.

\begin{lemma}\label{quasi-symmetric.lem}
If $\Lambda$ is quasi-symmetric \(resp.\ symmetric\) and the
action of $G$ on $\Lambda$ is trivial,
then $\Omega$ is quasi-symmetric \(resp.\ symmetric\).
\end{lemma}

\begin{proof}
  Taking the completion, we may assume that $R$ is complete.
  Then replacing $R$ by a Noether normalization of $R/\ann_R \Lambda$,
  we may assume that $R$ is a regular local ring, and
  $\Lambda$ is a faithful $R$-module.
  As the action of $G$ on $\Lambda$ is trivial, 
  $\Omega=\Lambda\otimes_R RG$ is quasi-symmetric (resp.\ symmetric),
  as can be seen easily.
\end{proof}

\paragraph
In particular, if $\Lambda$ is commutative quasi-Gorenstein
(resp.\ Gorenstein) and the action of $G$ on $\Lambda$ is
trivial, then
$\Omega=\Lambda G$ is quasi-symmetric (resp.\ symmetric).

\paragraph
In general, ${}_\Omega\Omega\cong \Lambda\otimes_R RG$ as $\Omega$-modules.

\begin{lemma}\label{isomorphic-Lambda-Omega.lem}
  Let $M$ and $N$ be right $\Omega$-modules, and let 
  $\varphi:M\rightarrow N$ be a homomorphism of right $\Lambda$-modules.
  Then $\psi:M\otimes RG\rightarrow N\otimes RG$ given by
  $\psi(m\otimes g)=g(\varphi(g^{-1}m))\otimes g$ is an $\Omega$-homomorphism.
  In particular,
  \begin{enumerate}
  \item[\bf 1] If $\varphi$ is a $\Lambda$-isomorphism, then
    $\psi$ is an $\Omega$-isomorphism.
  \item[\bf 2] If $\varphi$ is a split monomorphism in $\md \Lambda$,
    then $\psi$ is a split monomorphism in $\md\Omega$.
  \end{enumerate}
\end{lemma}

\begin{proof}
Straightforward.
\end{proof}

\begin{proposition}
  Let $G$ be a finite group acting on $\Lambda$.
  Set $\Omega:=\Lambda*G$.
  \begin{enumerate}
  \item[\bf 1] If the action of $G$ on $\Lambda$ is trivial and
    $\Lambda$ is quasi-symmetric \(resp.\ symmetric\), then so is $\Omega$.
  \item[\bf 2] If $\Lambda$ is pseudo-Frobenius \(resp.\ Frobenius\),
    then so is $\Omega$.
  \item[\bf 3] If $\Lambda$ is pseudo-quasi-Frobenius \(resp.\
    quasi-Frobenius\), then so is $\Omega$.
  \end{enumerate}
\end{proposition}

\begin{proof}
  {\bf 1} is Lemma~\ref{quasi-symmetric.lem}.
  To prove {\bf 2} and {\bf 3}, we may assume that $(R,\fm)$ is
  complete regular local and $\Lambda$ is a faithful module.

  {\bf 2}.
  \[
  (K_\Omega)_\Omega
  \cong
  \Hom_R(\Lambda\otimes_R RG,R)\cong \Hom_R(\Lambda,R)\otimes(RG)^*
  \cong
  K_\Lambda\otimes RG
  \]
  as right $\Omega$-modules.
  It is isomorphic to $\Lambda_\Omega\otimes RG\cong \Omega_\Omega$
  by Lemma~\ref{isomorphic-Lambda-Omega.lem}, {\bf 1}, since
  $K_\Lambda\cong \Lambda$ in $\md\Lambda$.
  Hence $\Omega$ is pseudo-Frobenius.
  If, in addition, $\Lambda$ is Cohen--Macaulay, then $\Omega$ is also
  Cohen--Macaulay, and hence $\Omega$ is Frobenius.

  {\bf 3} is proved similarly, using
  Lemma~\ref{isomorphic-Lambda-Omega.lem}, {\bf 2}.
\end{proof}

Note that the assertions for Frobenius and quasi-Frobenius properties also
follow easily from Lemma~\ref{SS.lem} and \cite[(3.2)]{SS}.

\section{Codimension-two argument}\label{codim-two.sec}

\paragraph Let $X$ be a locally Noetherian scheme, $U$ its open subscheme,
and $\Lambda$ a coherent $\O_X$-algebra.
Assume the $(S_2')$ condition on $\Lambda$.
Let $i:U\hookrightarrow X$ be the inclusion.
In what follows we use the notation for rings and modules to schemes and
coherent algebras and modules in an obvious manner.

\paragraph
Let $\M\in\md\Lambda$.
That is, $\M$ is a coherent right $\Lambda$-module.
Then by restriction, $i^*\M\in \md i^*\Lambda$.

\paragraph
  For a quasi-coherent $i^*\Lambda$-module $\N$, we have an action
  \[
  i_*\N\otimes_{\O_X}\Lambda
  \xrightarrow{u\otimes 1} i_*\N \otimes_{\O_X}i_*i^*\Lambda
  \rightarrow i_*(\N\otimes_{\O_U} i^*\Lambda)\xrightarrow{a} i_*\N.
  \]
  So we get a functor $i_*:\Mod i^*\Lambda \rightarrow \Mod \Lambda$, where
  $\Mod i^*\Lambda$ (resp.\ $\Mod \Lambda$) denote the category of quasi-coherent
  $i^*\Lambda$-modules (resp.\ $\Lambda$-modules).

\begin{lemma}\label{S_2-isom.lem}
  Let the notation be as above.
  Assume that $U$ is large in $X$ \(that is, $\codim_X(X\setminus U)\geq 2$\).
  If $\M\in (S_2')^{\Lambda,}$, then the canonical map
    $u:\M\rightarrow i_*i^*\M$ is an isomorphism.
\end{lemma}

\begin{proof}
  Follows immediately from \cite[(7.31)]{Hashimoto12}.
\end{proof}

\begin{proposition}\label{codim-two.prop}
  Let the notation be above, and let $U$ be large in $X$.
  Assume that there is a $2$-canonical right $\Lambda$-module.
  Then we have the following.
  \begin{enumerate}
  \item[\bf 1]
        If $\N\in (S_2')^{i^*\Lambda,U}$, then $i_*\N\in (S_2')^{\Lambda,X}$.
  \item[\bf 2] $i^*:(S_2')^{\Lambda,X}\rightarrow (S_2')^{i^*\Lambda,U}$ and
    $i_*:(S_2')^{i^*\Lambda,U}\rightarrow (S_2')^{\Lambda,X}$ are quasi-inverse
    each other.
  \end{enumerate}
\end{proposition}

\begin{proof}
  The question is local, and we may assume that $X$ is affine.

  {\bf 1}.
  There is a coherent subsheaf $\Q$ of $i_*\N$ such that $i^*\Q = i^*i_*\N=\N$
  by \cite[Exercise~II.5.15]{Hartshorne}.
  Let $\V$ be the $\Lambda$-submodule of $i_*\N$ generated by $\Q$.
  That is, the image of the composite
  \[
  \Q\otimes_{\O_X}\Lambda\rightarrow i_*\N\otimes_{\O_X}\Lambda\rightarrow i_*\N.
  \]
  Note that $\V$ is coherent, and $i^*\Q\subset i^*\V \subset i^*i_*\N=i^*\Q=\N$.

  Let $\C$ be a $2$-canonical right $\Lambda$-module.
  Let $?\da:=\uHom\Lop(?,\C)$, $\Gamma=\uEnd_{\Lambda}\C$,
  and
  $?\dda:=\uHom_{\Gamma}(?,\C)$.
  Let $\M$ be the double dual $\V\ddd$.
  Then $\M\in (S_2')^{\Lambda,X}$, and hence
  \[
  \M\cong i_*i^*\M\cong i_*i^*(\V\ddd)\cong i_*(i^*\V)\ddd\cong i_*(\N\ddd)\cong i_*\N.
  \]
  So $i_*\N\cong \M$ lies in $(S_2')^{\Lambda,X}$.

  {\bf 2} follows from {\bf 1} and Lemma~\ref{S_2-isom.lem} immediately.
\end{proof}

\end{document}